\RequirePackage{fix-cm}
\documentclass[twocolumn]{svjour3}          % twocolumn
\smartqed  % flush right qed marks, e.g. at end of proof
 \usepackage{mathptmx}      % use Times fonts if available on your TeX system
\usepackage{microtype}
\usepackage{natbib}
\usepackage{graphicx}
\usepackage{epsfig}
\usepackage{amssymb}
\usepackage{amsmath}
\usepackage{amsfonts}
\usepackage{subfig}
\usepackage{color}
\DeclareMathOperator*{\argmin}{arg\,min}
\journalname{arXiv.org}
\begin{document}
\title{On $1$--Laplacian Elliptic Equations Modeling Magnetic Resonance Image Rician Denoising\thanks{The first and second authors wish to thank the Spanish Ministerio de Econom\'ia y Competitividad for supporting Project TEC2012-39095-C03-02. The third author has partially been supported by the Spanish Ministerio de Econom\'ia y Competitividad under Project MTM2012-31103. Finally, the authors wish to thank the reviewers for providing helpful comments that have improved the final manuscript redaction.}}

\titlerunning{On $1$--Laplacian Elliptic Equations Modeling MRI Rician Denoising}        % if too long for running head

\author{\author{Adri\'an Mart\'{\i}n
		\and Emanuele Schiavi
		\and Sergio Segura de Le\'on}
}


\institute{
	A. Martin, E. Schiavi \at
	Departamento de Matem\'atica Aplicada, Ciencia y Tecnolog\'{\i}a de los Materiales y Tecnolog\'{\i}a Electr\'onica, Universidad Rey Juan Carlos, Tulip\'an S/N, M\'ostoles, Madrid, Spain (\email{adrian.martin@urjc.es, emanuele.schiavi@urjc.es})
           \and
           S. Segura de Le\'on \at
              Departament d'An\`alisi
              Matem\`atica, Universitat de Val\`encia,  Dr. Moliner 50, 46100
              Burjassot, Val\`encia, Spain (\email{sergio.segura@uv.es})
}

\date{July 2016}

\maketitle

\begin{abstract}
Modeling magnitude Magnetic Resonance Images (MRI) Rician denoising in a Bayesian or generalized Tikhonov framework using Total Variation (TV) leads naturally to the
consideration of nonlinear elliptic equations. These involve the
so called $1$--Laplacian operator and special care is needed to
properly formulate the problem. The Rician statistics of the data
are introduced through a singular equation with a reaction term
defined in terms of modified first order Bessel functions. An
existence theory is provided here together with other qualitative
properties of the solutions. Remarkably, each
positive global minimum of the associated
functional is one of such solutions. Moreover, we directly solve this non--smooth non--convex minimization problem using a convergent Proximal Point
Algorithm. Numerical results based on synthetic and real MRI demonstrate a better performance of the proposed method when compared to previous TV based models for Rician denoising which regularize or convexify the problem. Finally, an application on real Diffusion Tensor Images, a strongly affected by Rician noise MRI modality, is presented and discussed.
\keywords{$1$-Laplacian \and Total Variation operator \and Rician denoising  \and  Non-smooth non-convex energy minimization \and Global minimizer \and Magnetic Resonance Imaging \and Diffusion Tensor Imaging}
\end{abstract}

\section{Introduction}
Multiple applications in computer vision and digital image processing can be modeled from the field of quasilinear elliptic equations.
Variational formulations of these equations allow to introduce a concept of weak solution, which is well adapted to image analysis, providing faithful discontinuous solutions. Furthermore, the discrete formulations of these equations are readily suited for fast image processing.
In particular, medical image denoising is an important application which allows to reduce scanning time of the patients while preserving a good image quality. Moreover, several imaging applications like segmentation, classification, registration,
super-resolution, object recognition or tracking can benefit of pre--processed denoised images.

In this paper we focus on the modality of Magnetic Resonance
Imaging (MRI), where clinicians typically work with images contaminated by Rician noise. MRI scanners acquire complex data where both real and imaginary parts are corrupted with zero-mean uncorrelated Gaussian noise with equal variance. The calculation of the magnitude image transforms the original complex Gaussian noise into Rician noise \citep{Henkelman1985,Gudbjartsson1995}. The Rician distribution considerably differs from a Gaussian distribution when low signal-to-noise-ratio (SNR) data is considered. This is the reason why several denoising methods that take into account the Rician distribution of the noise are focused on Diffusion Weighted Images (DWI), one of the MRI modalities more severely affected by noise \citep{Basu2006,Wiest-Daessle2008a,Tristan-Vega2010}.

In particular, modeling these statistics in the
framework of a Tikhonov Regularization through the Total Variation (TV)
operator leads to consider a $1$--Laplacian elliptic equation
with a nonlinear lower order term defined in terms of modified
Bessel functions.

The TV operator
\begin{equation}\label{tv}
TV(u)=\displaystyle \mbox{sup}\left\{ \int_\Omega u\, \hbox{div\,}
\phi \, \Big|\, \phi \in C_c^{\infty} (\Omega ,\mathbb{R}^N
),  \| \phi \|_\infty \leq 1\right\}
\end{equation}	
was introduced in the image community by Rudin, Osher and Fatemi, \citep{Rudin1992a}
 through their celebrated denoising model (ROF in the following)
which is the Gaussian counterpart of the Rician model we are considering.
The $1$--Laplacian operator characterizes the subdifferential of  the TV functional; for a proof of such result in the $L^2$--framework, we refer to \citep{Andreu2004b} (see \citep[Proposition 1.10]{Andreu2004b}). We point out that in \citep{Bredies2012b} the $1$--Laplacian operator has also been characterized as the pointwise subdifferential of the TV operator
in form
\[ -\displaystyle  \mbox{div} \left( \frac{Du }{|Du |}\right)  \in \partial TV (u)\,. \]
It is well known that inverse ill--posed problems can be dealt
with in the framework of generalized Tikhonov regularization. The
resulting functional is composed of two basic terms which reflect
our belief in the data through one or more (hyper)--parameters
weighting the amount of regularization. This in turn determines
the smoothness of the denoised image and functional
analysis is invoked in order to select the appropriate functional
space. Sobolev spaces are rapidly ruled out because of their
excessive smoothing which generates continuous unrealistic images.
So this very nonlinear operator, the TV operator emerges because
it allows for (weak) distributional solutions in the very large
space of functions of bounded variation, those whose gradient is a
Radon measure \citep{Ambrosio2000}. Such a sophisticated setting is a generalized approach which allows for truly discontinuous
functions and opens the way to theoretical as well as practical
and accurate digital image processing \citep{Chambolle2010a}. Since the seminal paper from Rudin, Osher and Fatemi \citep{Rudin1992a}, there has been a burst in the application of the Total Variation
regularization to many different image processing problems which include inpainting, blind
deconvolution or multichannel image segmentation (see for instance \citep{Chan2005} for a review on this topic).
%In the last years, a generalized version of TV (the so-called Total Generalized Variation \citep{Bredies2010a}) which improves the original TV performance, alleviating the well known staircasing artifact introduced by TV processing, has grown as a widely
%used alternative in image processing applications.
Fast and robust numerical methods have been proposed to exactly solve convex optimization problems with TV regularization,
such as the dual approach of \citep{Chambolle2004} and, more
recently, the Split Bregman method \citep{Goldstein2009a} and the  primal--dual approach of \citep{Chambolle2011}.

Our proposed model equation arises as the
(formal) Euler--Lagrange equation associated to an energy
minimization problem obtained in a Bayesian framework.
A key feature of this problem is that the nonlinear term modeling Rician noise  in the energy functional
can be a non--convex  changing sign function  with a double well profile. This leads to the study of non--convex non--smooth minimization problems.
In fact, this non--convexity of the energy functional is crucial because otherwise we could show uniqueness of the trivial solution $u\equiv 0$.
The variational minimization problem associated to the model equation we consider
was proposed in \citep{Martin2011}, where the multivalued Euler--Lagrange equation for the $1$--Laplacian operator is deduced as a first order
necessary optimality condition. This minimization problem was simultaneously and independently considered in \citep{Getreuer2011}, where blurring effects were included
and existence and comparison results in the pure denoising case were reported. In order to cope with the multivalued Euler--Lagrange equation an $\epsilon$-regularization of the TV term was introduced in both works \citep{Martin2011,Getreuer2011} .
More recently, in \citep{Chen2015}
a convex variational model for restoring blurred images corrupted by Rician noise have been proposed to overcome the difficulties related to the non--convex nature
of the original problem we are considering here.

The non--smoothness property of the model comes from
the very singular $1$--Laplacian elliptic equation, which had firstly been studied as a limit of equations
involving the $p$--Laplacian.
The interest in studying such a case
came from an optimal design problem in the theory of torsion and
related geometrical problems (see \citep{Kawohl1990} and \citep{Kawohl1991} for
constant data, and  \citep{Cicalese2003} for more general data). The suitable notion of solution to the $1$--Laplacian had to wait at the turn
of the century \citep{Andreu2001}. Other important related papers
published in the early twenty--first century include
\citep{Andreu2001a,Bellettini2002,Demengel2002,Kawohl2003,Andreu2004, Demengel2004,Bellettini2006}.
Due to its unique properties, this operator has been the
optimal choice for PDE based image processing in the last twenty
years. Briefly, the $1$--Laplacian describes isotropic diffusion
within each level surface with no diffusion across different level
surfaces. In this way, its action does not over-regularize the data and preserves
edges and fine details. This is not true when the $p$--Laplacian operator, for $p>1$, is used, since an artificial smoothing is introduced.

While the ROF model has been mathematically studied
and existence and uniqueness results have been obtained
\citep{Chambolle1997},  the Euler--Lagrange quasilinear equation associated to the Rician problem has not been considered yet for mathematical analysis.
Notice that the same is true even for the semilinear equation accounting for Rician noise and
linear diffusion. Here we focus on the mathematical analysis of the TV based Rician model. We show that the nonlinear $1$--Laplacian problem has, aside from the trivial solution, at least a positive distributional solution
which is also a global minimum of the energy problem (provided that the datum is big enough). This result makes the solutions of the TV Rician denoising model attractive for the application in MRI and in particular for the DWI modality.
%and practically computable as we shall see in the numerics.
The existence result is based on the consideration of a sequence
of approximating problems of the $p$--Laplacian type for which no
existence results are known due to the very special nonlinearity
associated to the Rician noise term. Standard techniques can be
used. When $p=2$, existence and uniqueness of positive solutions is
also deduced. In contrast, for general $1\leq p <2$, the uniqueness of positive
solutions is still an open problem. Nevertheless, it is proved that for constant data we
have uniqueness of constant solutions for any $p$.

The numerical resolution of the proposed model is also challenging
because the energy functional is non--convex for any $p$ and also
non--smooth for $p=1$. To cope with the non--convexity we propose
a suitable decomposition of the energy functional, which allows to
write it as a Difference of Convex (DC) functionals. A
primal--dual approach (suitable for non differentiable energy
functionals such as the TV operator) embedded into a proximal
algorithm (suitable for DC functionals) is then applied to show,
also numerically, the convergence of the $p$--Laplacian
approximate solutions to the true $1$--Laplacian solution when $p
\to 1$. This provides a unified framework in which these problems
can be solved using the same algorithm and then fairly
compared. Our numerical method is then successfully compared with the primal gradient descent algorithm presented in \citep{Getreuer2011} and the convexified models of \citep{Getreuer2011} and \citep{Chen2015}.

This paper is organized as follows. In Section 2 we define the
model problem characterizing the Bessel ratio function and its
properties jointly with the statement of our main result (see
Theorem \ref{main} below). Weak solutions are defined in Section 3
where the main result is obtained considering suitable
regularizing approximating problems of the $p$--Laplacian type.
Some qualitative properties are discussed in Section 4, before the numerical resolution of the related minimization problem is presented in Section 5. Finally, in section 6, the performance of the algorithm is compared to other related methods and an application on real DTI is presented.

\section{Preliminaries}

\subsection{The model problem and the statement of the main result}%\label{s1}
Let $\Omega$ be an open, bounded domain in $\mathbb{R}^N$
($N\ge2$) with Lipschitz boundary $\partial\Omega$ (usually a
rectangle in image processing). Thus, there exists a outer unit
normal vector $n(x)$ at $x\in\partial \Omega$, for $\mathcal
H^{N-1}$--almost all point; here and in what follows $\mathcal
H^{N-1}$ stands for the $(N-1)$--dimensional Hausdorff measure.

We will consider in $\Omega$ a Neumann problem involving the
$1$--Laplacian. This operator has to be studied in the framework
of functions of bounded variation. Recall that a function
$u\>:\>\Omega\to\mathbb R$ is said to be of bounded variation if
$u\in L^1(\Omega)$ and its distributional gradient $Du$ is a
(vector) Radon measure having finite total variation. We denote by
$BV(\Omega)$ the space containing all functions of bounded
variation. For a systematic study of this space, we refer to
\citep{Ambrosio2000} (see also \citep{Giusti1984}). The appropriate concept of
solution to deal with the Neumann problem for the $1$--Laplacian
is introduced in \citep{Andreu2001}. For a review on the early
development of variational models in image processing and a deep
study of equations involving the $1$--Laplacian, see \citep{Andreu2004}.

The boundary value problem in which we are interested is:
\begin{equation}\label{hbes}
\left\{
\begin{array}{ll}
-\displaystyle  \mbox{div} \left( \frac{Du }{|Du |}\right)+h' (x, u)=0,&\mbox{in}\,\Omega , \\[0.5cm]
\displaystyle \left( \frac{Du }{|Du |}\right) \cdot n =0, &
\mbox{on}\, \partial \Omega\,.
\end{array}
\right.
\end{equation}
We shall assume that $h':\, \Omega \times \mathbb{R}\to
\mathbb{R}$ is a non monotone Carath\'eodory function defined as
\begin{equation}\label{semi}
h' (x,u)=\left(\displaystyle \frac{\lambda}{\sigma^2} \right) u -\left(\displaystyle \frac{\lambda}{\sigma^2} \right)  r_\sigma (x,u) f(x)
\end{equation}
%$$
%\left(\displaystyle \frac{\lambda}{\sigma^2} \right) u- \displaystyle
%\left( \frac{\lambda}{\sigma^2}\right)\left[  \frac{I_1 (uf /\sigma^2 )}{I_0 (uf /\sigma^2 )} \right] f(x)
%$$
where $\lambda >0$ and $\sigma^2 \neq 0$ are real given parameters, $f(x)\ge 0$ for almost all $x\in \Omega$, and
the function
\begin{equation}\label{ratio}
\displaystyle
r_\sigma (x,u)=\displaystyle \frac{I_1 \displaystyle \left( \frac{u (x) f(x) }{\sigma^2 }\right) }{I_0 \displaystyle \left( \frac{u (x)f(x)}{ \sigma^2 } \right)}
\end{equation}
is the ratio between the first and zero order modified  Bessel
functions of the first kind. Series representations and general properties can be found in \citep{Watson1922}.
Notice the dependence (that we shall omit) $r_{\sigma} (x,u)=r(x,u)$ of the Bessel ratio function on the parameter $\sigma^2$, which is the estimated variance of the original Gaussian noise of the complex MRI data. This implicit dependence renders problem (\ref{hbes}) a truly $2$--parametric problem in so far $\sigma^2$ cannot be scaled out from $\lambda$ and it has to be estimated
from the noisy data $f(x)$.
%This semilinearity act as a reaction term in the equation.

Assuming $\lambda >0$, $\sigma^2 \neq 0$ fixed and $f(x) \in L^{\infty} (\Omega )$ given, problem (\ref{hbes}) reads:
\begin{equation}\label{bes}
\left\{
\begin{array}{ll}
\displaystyle \left( \frac{ \lambda}{\sigma^2}\right) u  -\displaystyle \mbox{div} \left( \frac{Du }{|Du |}\right)=\left(\displaystyle \frac{ \lambda}{\sigma^2} \right) r(x,u) f ,&\mbox{in }\,\Omega\,; \\[0.5cm]
\displaystyle \left( \frac{Du }{|Du |}\right) \cdot n =0, & \mbox{on}\, \partial \Omega .
\end{array}
\right.
\end{equation}
The modified Bessel functions $I_\nu (s)$, $\nu \geq 0$, $s\geq 0$ which define the ratio $r(x,u)$ \eqref{ratio}
are monotone, exponentially growing functions and this distinguish
their behavior from ordinary Bessel functions which have
oscillating wave--like forms \citep{Amos1974,Neuman1992}. Moreover $I_0 (0)=1$, $I_0 (s)
>1$ for any $s>0$ and  $I_\nu (0)=0$,  $I_\nu (s)>0$ for any $s>0$
and $\nu
>0$ so $r(x,0)=0$ and the Bessel ratio function $r (x,u)$ in (\ref{ratio}) is well--defined and non--negative for
any $f\geq 0$ and $u \geq 0$. Also $I_1 (s)<I_0 (s)$ for any $s>0$ and then $0\leq r(x,u) <1$.
By (\ref{semi}) we then have $h' (x,0)=0$ and
$u\equiv 0$ is always a solution of (\ref{hbes}) and (\ref{bes})
for any non-negative datum $f (x)$ and fixed parameters $\lambda >0$ and $\sigma^2 \neq 0$.

The specific form of $h' (x,u)$ given in \eqref{semi} describes the Rician noise distribution of a given datum image $f(x)$ and it has been deduced in several papers
dealing with medical imaging since the paper \citep{Basu2006} where it was proposed for DTI.
When dealing with the image processing application we shall assume that $f\in L^\infty (\Omega )$ even if
our existence theory applies more generally to $f\in L^2 (\Omega )$.

The function $h' (x,u)$ is the Gateaux derivative of $h(x,u)$:
\[ h(x,u)=\displaystyle \int_0^u h' (x,t )dt \]
Using (\ref{semi}) we have:
\begin{equation}\label{hpr}
h(x,u)=\displaystyle \left( \frac{\lambda }{2\sigma^2}\right) u^2 -\lambda \log I_0 \left( \frac{u f}{\sigma^2}
\right)
\end{equation}
with $h(x,0) =0$ and the logarithm is well--defined and
nonnegative for any $f\geq 0$ and $u \geq 0$ because of $I_0
(s)\geq 1$, $\forall\,s\geq 0$.

Following the Bayesian modeling approach, the associated
minimization problem is
\begin{equation*}%\label{minp}
\min_{u \in BV(\Omega)}
J_1 (u)+ H(u,f)\,,
\end{equation*}
where $J_1 (u)=TV(u)$ is the Total Variation regularization
functional, previously defined in \eqref{tv}, and that can also be denoted as
\begin{equation*} %\label{JU}
J_1 (u)=\int_\Omega |D u|\,.
\end{equation*}
The fidelity term (modelling Rician noise) is
\begin{eqnarray}\label{hener}
H(u,f)&=&\displaystyle \int_\Omega h(x,u)dx \nonumber\\
&=&\displaystyle \displaystyle \left( \frac{\lambda }{2\sigma^2}\right) \int_\Omega u^2 dx -
\lambda \displaystyle \int_\Omega \log I_0 \left( \frac{u f}{\sigma^2} \right)dx\,.
\end{eqnarray}
Notice that $H(0,f)=0$, $\forall\,f$.

The minimization problem for image denoising of Rician corrupted
data is formulated as follows.  An equivalent formulation is
considered in \citep{Getreuer2011}.  Fixed real parameters $\lambda >0$ and
$\sigma^2 \neq 0$ and given a noisy image $f \in L^{\infty}
(\Omega)$ recover a clean image $u \in BV(\Omega)\cap L^{\infty}
(\Omega)$ minimizing the energy:
\begin{eqnarray}\label{minFun}
E_1 (u)&=&J_1 (u)+ H(u,f) \\
&=&\int_\Omega  |D u|  + \lambda  \int_\Omega \frac{u^2}{2\sigma^2} -\lambda\log I_0 \left( \frac{u f}{\sigma^2} \right) dx \,.\nonumber
\end{eqnarray}
This minimization problem can naturally be studied in the $L^2$--setting since
\[
|H(u,f)|\le \frac{\lambda}{2\sigma^2}\int_\Omega u^2\,dx+\frac{\lambda}{\sigma^2}\int_\Omega |u|f\,dx
\]
(see \eqref{quad} in Lemma \ref{prop-h} below).
Thus, our main result can be stated
as follows:

\begin{theorem}\label{main}
	Let $\lambda >0$ and $\sigma^2 \neq 0$ be given real parameters. For
	every non--negative $f\in L^2(\Omega)$, there exists a
	non--negative $u\in BV(\Omega)\cap L^2(\Omega)$ which is a
	solution to problem \eqref{hbes},
{\color{black}
in the sense of Subsection 3.1,
}
 and it is a {\em global} minimum
	of functional $E_1$ in (\ref{minFun}).
\end{theorem}

\begin{remark}
This existence result relates problem \eqref{hbes}  and the
	global minimization of functional \eqref{minFun}, which is a
	non--smooth and non--convex optimization problem.
  Its proof can be found in Section 3 below, while Section 4 is devoted to complete this theorem. Among others, it is shown that the solution we find satisfies the following properties:
\begin{enumerate}
  \item If $f\in L^\infty(\Omega)$, then $u\in L^\infty(\Omega)$ and $\|u\|_\infty\le\|f\|_\infty$.
    \item Solution $u$ vanishes identically when $f(x)\le\sqrt{2\sigma^2}$ a.e. in $\Omega$.
  \item Solution $u$ is strictly positive when $f(x)\ge\mu>\sqrt{2
		\sigma^2} $ a.e. $x \in \Omega$, and moreover $E_1(u)<0$ holds.
\end{enumerate}
This last feature provides a sufficient condition in order to have a non
	trivial minimizer of functional \eqref{minFun}.
\end{remark}

{\color{black}
\begin{remark}
The problem of minimizing $E_1$ has also been considered by Getreuer et al. (their results were announced in \citep{Getreuer2011} and proved in \citep{Getreuer2011a,Tong2012}). It is worth comparing these results with those in the present paper since
 both approaches are very different. We prove our results through the formal Euler--Lagrange equation of the minimization problem,
while the results in \citep{Getreuer2011a,Tong2012} are obtained by direct methods. We explicitly point out two aspects:
\begin{enumerate}
  \item Our existence result is more general, since we take data $f$ belonging to $L^2(\Omega)$,
  and \citep{Getreuer2011a,Tong2012} consider data $f\in L^\infty(\Omega)$ with the additional assumption
  $\inf_{x\in\Omega} f(x)\ge\alpha>0$.
  \item One important feature of the present paper is that a simple condition is provided to distinguish data which lead to non trivial solutions.
  Instead, the results by Getreuer et al. do not identify non trivial solutions.
\end{enumerate}
\end{remark}
}

\begin{remark}\label{rm21}
{\color{black} Uniquenes of non-trivial solutions is still an open problem.}
Using the same
	arguments from \citep{Aubert2008} and some properties of the modified
	Bessel Functions, a comparison result for the solutions of the
	minimization problem is stated in \citep[Theorem 2]{Getreuer2011}
{\color{black}
and proved in \citep{Getreuer2011a,Tong2012}.
}
 This comparison
	result establishes that given
	$0<f_1 < f_2$ a.e. $x\in \Omega$, then $u_1 \leq u_2$ a.e. $x\in \Omega$, with $u_1$, $u_2$ being minimizers of \eqref{minFun}  for $f=f_1$, $f=f_2$ respectively.
	Since $f_1$ and $f_2$ must be different, it does not imply
	uniqueness. {\color{black} 
	Some partial results about uniqueness of non trivial solutions shall be presented in section 4.2}
\end{remark}

The existence result in Theorem \ref{main} will be proved by
approximating our functional through functionals defined on the
Sobolev space $W^{1,p}(\Omega)$ and having $p$--growth (with
$p>1$). The main advantage of these approximating functionals is
their differentiability (in contrast with $E_1$, which is not
differentiable). So, we introduce, for subsequent analysis, the
(differentiable) energy
\begin{eqnarray}\label{ep}
E_p (u) &=& J_p (u)+ H(u,f) \\
&=&\displaystyle \frac1p \int_\Omega
|\nabla u |^{p}\, dx +\lambda \int_\Omega \frac{u^2}{2\sigma^2} - \lambda\log I_0 \left( \frac{u f}{\sigma^2} \right) dx \,.\nonumber
\end{eqnarray}
Notice that $E_p (0)=0$ for any $p> 1$, and also $E_1 (0)=0$. The weak (distributional)
solutions of \eqref{bes} formally coincide with the critical
points of \eqref{minFun}.  The crucial point is that these
energies (including  \eqref{ep} for $p>1$) may be non--convex
depending on the datum $f$ and the (estimated) parameter $\sigma^2$.
This fact does not depend on the regularizer but it is a feature of the
Rician likelihood function. To explore further this point we
analyse the behaviour of the Bessel ratio function defined in
\eqref{ratio} which governs the qualitative properties of the energies \eqref{minFun} and \eqref{ep}.
This leads to show the coercitiveness of the functional in section \ref{subsect_coerc},
implying the existence result for the $p$-approximating problems in section \ref{subsec:approx}.
\subsection{A Non--Convex Semi--Linearity}
The characterization of the model semilinearity $h(x,u)$ leads to
the study of the properties of the modified Bessel functions of
the first kind. Our results are founded on some fundamental
inequalities regarding the ratio function $r(x,u)$ and its
derivative which can be found in \citep{Amos1974}.
These results will allow to characterize suitable growth conditions
related to the Rician statistics. Moreover we shall prove that, depending on the data and parameters of the problem,  $h^{''} (x,u)$ is negative near $u=0$ and hence $h'$ is non--monotone
and $h$ is non--convex.

%It then models a sublinear reaction--absorption balance for
%$1-$Laplacian diffusion.

\begin{lemma}\label{prop-h} Let $h'$ be defined as
	in (\ref{semi}) with datum $f(x)\geq 0$ and fixed parameters $\lambda >0$ and $\sigma^2 \neq 0$. Then
	\begin{equation}\label{quad}
	|h (x,u)| \leq \displaystyle \frac{\lambda}{2\sigma^2}  u^2+\frac{\lambda}{\sigma^2}|u|f(x),\quad a.e.\,x\in\Omega .
	\end{equation}
	\begin{equation}\label{sub}
	|h' (x,u)| \leq \left(\displaystyle \frac{\lambda}{\sigma^2} \right) (| u |+f(x)),\quad a.e.\,x\in\Omega .
	\end{equation}
	and
	\begin{equation}\label{noconv}
	\left(\displaystyle \frac{\lambda}{\sigma^2} \right) \left( 1 - \displaystyle \frac{f^2 }{2\sigma^2 } \right)
	\leq h^{''}(x,u) \leq \left(\displaystyle \frac{\lambda}{\sigma^2} \right),\quad a.e.\,x\in\Omega .
	\end{equation}
	Moreover	
	\begin{equation}\label{h2}
	\displaystyle \lim_{u\to 0^+} h^{''}(x,u) = h^{''}(x,0) = \left(\displaystyle \frac{\lambda}{\sigma^2} \right) \left( 1 - \displaystyle \frac{f^2 }{2\sigma^2 } \right),
	\end{equation}
	and
	\begin{equation}\label{h2inf}
	%\displaystyle \lim_{u\to 0^+} h^{''}(x,u) = h^{''}(x,0) = \left(\displaystyle \frac{\lambda}{\sigma^2} \right) \left( 1 - \displaystyle \frac{f^2 }{2\sigma^2 } \right),\quad
	\displaystyle \lim_{u\to \infty} h^{''}(x,u) =\displaystyle
	\frac{\lambda}{\sigma^2} \,,\quad a.e.\,x\in\Omega .
	\end{equation}
\end{lemma}

{\sc Proof:} {\color{black} In order to simplify the notation when using the results of \citep{Amos1974}, we define $s= u (x)f(x)/\sigma^2 $ for fixed $x\in\Omega $ and denote the ratio function  $r(x,u)=r(s)=I_1 (s)/I_0 (s)$, $s\geq 0$. Please notice that $r'(x,u )$ is the Gateaux derivative while $r' (s)$ is the derivative w.r.t
the real, non-negative parameter $s$.}

By definition
and the monotonicity properties of the modified Bessel functions
$0\leq r(s)<1$ for any $s>0$ and $r(s)\to 1$ when $s\to \infty$.
The first inequality  is then straightforward.   We simply use
definition (\ref{semi}), the fact that $0\leq r(x,u) <1$ and the
triangle inequality to deduce that $h'$ verifies the sublinear
growth condition \eqref{sub} for a.e $x\in \Omega$.

As a consequence, we obtain \eqref{quad}. Indeed,
\begin{multline*}
  |h (x,u)| \le \int_0^{|u|}|h'(x,t)|\, dt\\
  \le \frac{\lambda}{\sigma^2} \int_0^{|u|} (| t |+f(x))\, dt
  \le \frac{\lambda}{\sigma^2}\Big( \frac{u^2}2+|u| f(x)\Big)\,.
\end{multline*}
%Notice that
%$$
%-\left(\displaystyle \frac{\lambda}{\sigma^2} \right) r (x,u) f(x)\leq  h' (x,u) \leq
%\left(\displaystyle \frac{\lambda}{\sigma^2} \right)  u
%$$
In order to show (\ref{noconv}) we compute
%\vspace{2.0cm}
the second derivative of $h(x,u)$ with respect to $u$ which reads
\begin{eqnarray}\label{ddh}
h^{''} (x,u)&=&\left(\displaystyle \frac{\lambda}{\sigma^2} \right)
[1-  \displaystyle r'(x,u) f(x)] \\
&=&\frac{\lambda }{\sigma^2}  \displaystyle \left[1- \left(\frac{f^2
}{2\sigma^2}\right) \left( 1+ \frac{ I_2 \displaystyle \left(
\frac{u f}{\sigma^2} \right)}{ I_0 \displaystyle \left( \frac{u
	f}{\sigma^2} \right)}          -2  \frac{ I_1^2 \displaystyle
\left( \frac{u f}{\sigma^2} \right)}{ I_0^2 \displaystyle \left(
\frac{u f}{\sigma^2} \right)}  \right)\right]\nonumber
\end{eqnarray}
where we used that $I_0^{''} (s)=I_1^{'}
(s)=(1/2)[I_2 (s)+I_0 (s) ]$, $\forall\,s\geq 0$ \citep{Amos1974}.
%\begin{equation}\label{local}
%h^{''} (x,0)=%\frac{\lambda }{\sigma^2}  -\displaystyle \frac{\lambda }{2}\left( \frac{f}{\sigma^2}\right)^2=
%\frac{\lambda }{\sigma^2}  -\displaystyle \frac{\lambda }{\sigma^2}\left( \frac{f^2 }{2\sigma^2}\right) =
%\displaystyle \left( \frac{\lambda }{2\sigma^4}\right) ( 2\sigma^2 -f^2  )
%\displaystyle \left( \frac{\lambda }{\sigma^2} \right)\left( 1-\frac{f^2 }{2\sigma^2}\right)
%\end{equation}

Reasoning as in \citep{Amos1974} and using its formulas 11, 12, 15, pg.242, the following bounds hold:
\begin{equation*}%\label{amos1}
0\leq \displaystyle \frac{s}{1+\sqrt{s^2+1}} \leq r(s)\leq
\displaystyle \frac{s}{\sqrt{s^2+4}} <1 ,\quad s\geq 0.
\end{equation*}
Using that
$\displaystyle r'(s)=1-\displaystyle \frac{r(s)}{s}-r^2 (s)$, $s>0$,
and inequalities $0< \displaystyle r' (s) < \displaystyle \frac{r(s)}{s}$, $s>0$ (formula 15 in \citep{Amos1974})  we get the improved bounds
\begin{equation}\label{amos3}
0 < r' (s) < \displaystyle \frac{1}{1+\sqrt{s^2+1}} <
\frac{r(s)}{s} < \frac{1}{\sqrt{s^2+4}}<\frac{1}{2}\,,
\end{equation}
for all $s>0$.

To show \eqref{noconv} we derive with respect
to $u$ the relationship $r(x,u)=r(s)$ to have
\[f(x)r' (s) =\sigma^2 r'(x,u)\]
and
%0 \leq \sigma^2 \frac{\partial r}{\partial u} (x,u) \leq \displaystyle f (x)\frac{r(s)}{s} =\displaystyle \sigma^2 \frac{r(x,u)}{u} <f (x) ,\quad a.e.\, x\in \Omega
%0 \leq  \frac{\partial r}{\partial u} (x,u) \leq \displaystyle \frac{r(x,u)}{u} <\frac{f (x)}{\sigma^2} ,\quad a.e.\, x\in \Omega
\[ 0 \leq \sigma^2 r' (x,u) =  f (x)r' (s) < \frac{f(x)}{2},\quad \mbox{a.e. }\, x\in \Omega \]
%and $r(x,u)$ is not Lipschitz (for $\sigma^2 \to 0$ the bound $f/\sigma^2 $ can blow up).
because $\displaystyle r'(s) < \frac12$ for any $s\geq 0$ by \eqref{amos3}. Using
(\ref{ddh}) the above inequality implies:
%
%$$
%\left(\displaystyle \frac{\lambda}{\sigma^2} \right) \geq \left(\displaystyle \frac{\lambda}{\sigma^2} \right)  -\left(\displaystyle \frac{\lambda}{\sigma^2} \right)  \frac{\partial r}{\partial u} (x,u) f(x) \geq
%\left(\displaystyle \frac{\lambda}{\sigma^2} \right) - \displaystyle \left(\displaystyle \frac{\lambda}{\sigma^2} \right)f(x)\frac{r(x,u)}{u} >
%$$
%$$
%> \left(\displaystyle \frac{\lambda}{\sigma^2} \right) - \displaystyle \left(\displaystyle \frac{\lambda}{\sigma^2} \right)\frac{f^2 (x)}{\sigma^2 } =
% \left(\displaystyle \frac{\lambda}{\sigma^2} \right) \left[ 1 - \displaystyle \left(\displaystyle \frac{f^2 }{\sigma^2 } \right)\right]
%$$
%and $h^{''} (x,u)$ can be negative so $h' (x,u)$ is not necessarly monotone. Multiple solutions may exist.
%which is
%
%We consider (\ref{ddh})
%
%
\begin{align*}
\left(\displaystyle \frac{\lambda}{\sigma^2} \right) \left[ 1 -
\displaystyle \left(\displaystyle \frac{f^2 }{2\sigma^2 }
\right)\right] &\leq h^{''} (x,u) =\\
&=\left(\displaystyle \frac{\lambda}{\sigma^2}
\right) \left[ 1 - r'(x,u) f(x)
\right]  \le \displaystyle \frac{\lambda}{\sigma^2}
\end{align*}
and (\ref{noconv}) holds true. Finally (\ref{h2}) is checked using (\ref{ddh})
and $I_0 (0)=1$, $I_\nu (0)=0$ for any $\nu >0$.  Because of  $I_2 (s)/I_0 (s) \to 1$, $I_1 (s)/I_0 (s) \to 1$ when $s\to \infty$ we deduce (\ref{h2inf}).
{\quad\rule{2mm}{2mm}\medskip}

\begin{figure}[h!tb]\label{figh2}
	\begin{center}
		{\includegraphics[width=0.9\columnwidth]{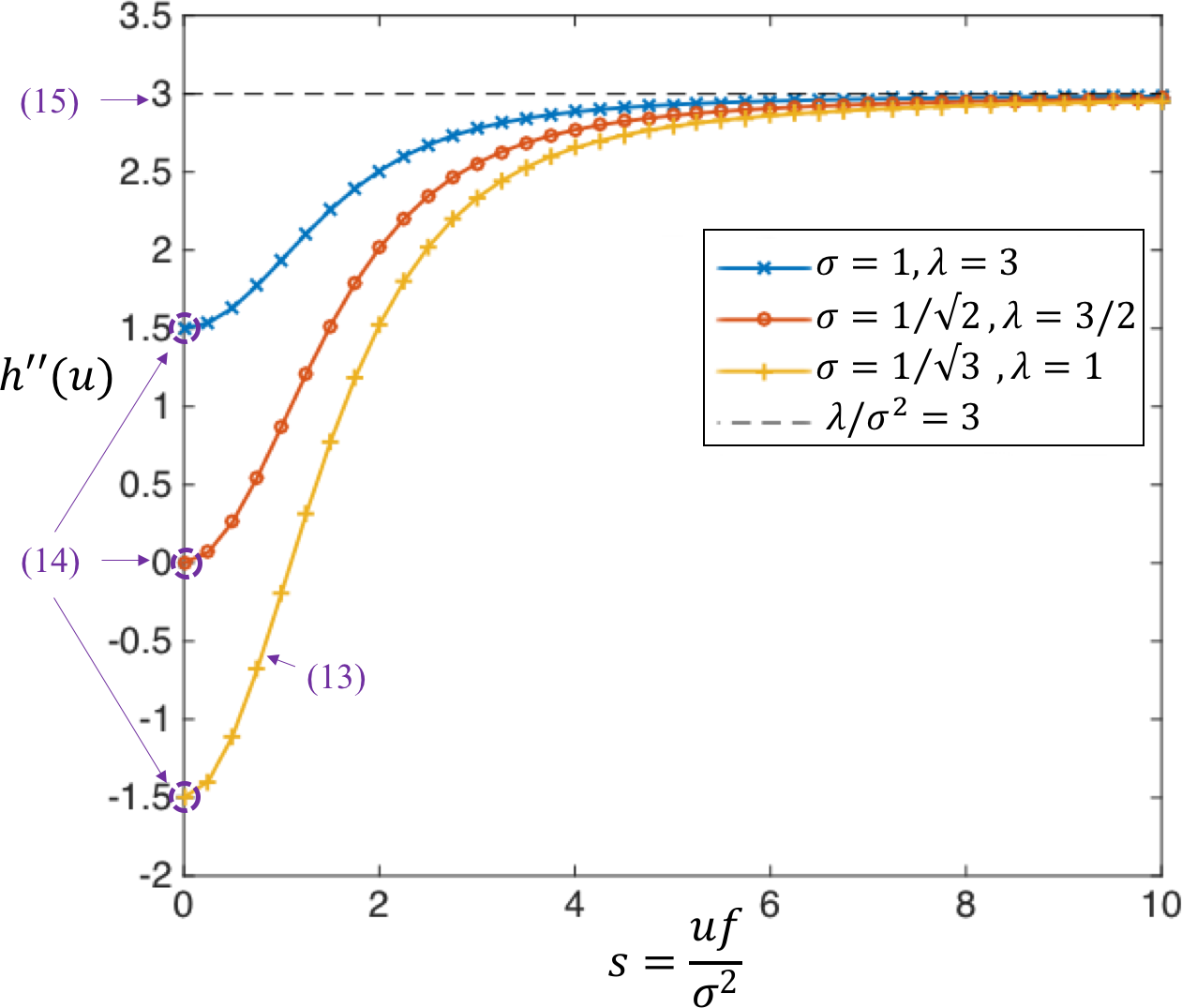}}
	\end{center}
	\caption{The profile of $h^{''}(x,u)$ is computed for constant data $f=1$ and
		parametric  values $f^2 /\sigma^2=1$, $f^2 /\sigma^2=2$ and $f^2 /\sigma^2=3$.
		The values of $\lambda$ is chosen to get a constant ratio $\lambda/\sigma^2=3$. A limit behavior is obtained when $f^2 /\sigma^2=2$ ($\sigma=1/\sqrt{2}$, in red).
		For $f^2 /\sigma^2\leq 2$ we have uniqueness. On the other hand, for $f^2 /\sigma^2>2$  we have $f^2
		>2\sigma^2$ and the corresponding profile is negative in a
		neighborhood of $s=0$. Some properties of $h''(x,u)$, \eqref{noconv}, \eqref{h2} and \eqref{h2inf}, can be observed in the figure.}
\end{figure}

\begin{figure}[h!]
	\begin{center}
		{\includegraphics[width=0.97\columnwidth]{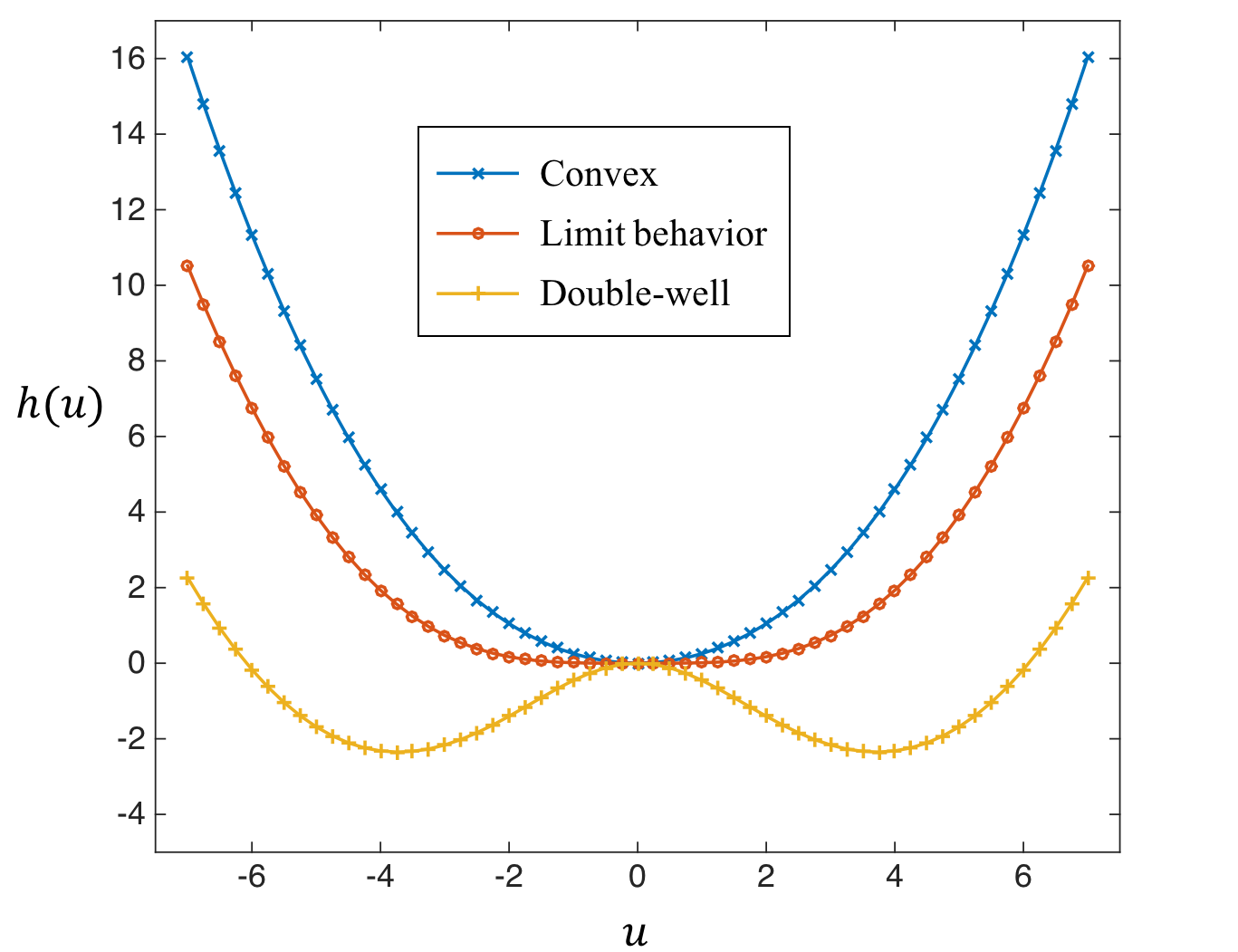}}
	\end{center}
	\caption{The double--well potential for parametric  values
		$\lambda =\sigma^2 =5$ is obtained when $f^2 >2\sigma^2=10$. In
		the figure above we represent the profile of function $t\mapsto
		h(t)$ for
		$f^2=\sigma^2 =5$ (Convex case), $f^2=2\sigma^2 =10$ (limiting behavior),
		$f^2=4\sigma^2 =20$ (double well).}
	\label{fig1}
\end{figure}

As a consequence of the above analysis $h^{''}(x,u)>0$ a.e. in $\Omega$ for (uniformly) small data
$f(x) <\sqrt{2\sigma^2}$ and then $h'$ is monotone increasing and uniqueness
of the trivial solution can be deduced (see Subsection \ref{uniqu}
below).
%The interesting case is when
%$\Omega_{\sqrt{2\sigma^2}}
%=\{ x/\, f(x)
%>\sqrt{2\sigma^2}\}$ has positive measure, which may easily be
%verified by checking $\|f\|_2^2>2\sigma^2|\Omega|$. In this case
%$h'$ is non--monotone and non--trivial solutions can exist.
Summing up we have shown (see figure \ref{fig1}) that $h^{''} (x,0) >0$ for $f (x) <\sqrt{2\sigma^2} $, $h^{''} (x,0) =0$ for $f (x) =\sqrt{2\sigma^2} $ and $h^{''} (x,0) <0$ for $f (x) >\sqrt{2\sigma^2} $.
% $a.e$ in $\Omega $.
The same fact is true for small $u$ as $h''(x,u)$ is continuous
with respect to $u$. These properties characterize the local
behavior near $u=0$ of $h(x,u)$. It turns out that $h^{''} (x,u)$
is a changing sign function depending on the datum $f$ and the
parameter $\sigma^2$. Then $h(x,u)$ is possibly non-convex.
% depending on the data $f$
%and $\sigma^2 $ (its noise variance).
This implies that $h' (x,u)$ is non monotone. Multiple solutions
to problem \eqref{hbes} corresponding to critical points of the
energy functional may exist. For $f\equiv 0$ we have
$h(u)=(\lambda /2\sigma^2 ) u^2 $, $h' (u)=(\lambda /\sigma^2 )u$
and $h^{''} (u) =\lambda /\sigma^2
>0$. Multiplying (formally) by $u$ in the model equation appearing
in (\ref{hbes}) and integrating it is easily seen that $u\equiv 0$
is the unique solution. We shall see in Subsection \ref{uniqu} that
the same phenomenon is true when $f$ is small enough.

%\subsection{Profile of function $h$}
To get a deep insight into the features of the energy term related
to Rician noisy data, in this subsection, we fix $x\in\Omega$ and
describe the profile of $h(x,u)$ defined in (\ref{hpr}). We have:

\begin{lemma}\label{lemma2}
	Let $h$ be defined as
	in (\ref{hpr}) with datum $f(x)\geq 0$, a.e. $x\in\Omega $ and
	fixed parameters $\lambda >0$ and $\sigma^2 \neq 0$. Then
	\begin{equation}\label{linfty}
	\lim_{t\to\pm\infty}h(x,t)=+\infty , \quad a.e.\,x\in\Omega .
	\end{equation}
	Moreover:
	\begin{enumerate}
		\item If $f(x)^2\le 2\sigma^2$ a.e.
		$x\in\Omega $, then the function $t\mapsto h(x,t)$ is convex and
		its minimum is attained at $0$.
		\item If $f(x)^2> 2\sigma^2$ a.e.
		$x\in\Omega $, then $t\mapsto h(x,t)$ has a unique positive
		critical point where it attains a global minimum.
	\end{enumerate}
\end{lemma}

{\sc Proof:} We fix $x\in\Omega $. When $f(x)=0$, the result is straightforward since then $h(x,t)=\left(\displaystyle \frac{\lambda}{2\sigma^2} \right)t^2$.

Assuming that $f(x)>0$, we begin by showing the limit behavior.  Consider
\eqref{semi} written in form
\[ h' (x,u)+\left(\displaystyle \frac{\lambda}{\sigma^2} \right)
r(x,u) f(x)= \left(\displaystyle \frac{\lambda}{\sigma^2} \right)
u\,. \]
As $h' (x,u)$ is the Gateaux derivative of $h (x,u)$ we formally
integrate in $(0,|u|)$ with respect to $u$ to obtain
\[ h (x, |u|)+\left(\displaystyle \frac{\lambda}{\sigma^2} \right)
\int_0^{|u|} r(x,t) f(x)dt = \left(\displaystyle
\frac{\lambda}{2\sigma^2} \right) |u|^2 \,. \]
We deduce from the boundedness $| r (x,t) | \leq 1$ a.e. in $\Omega $
for any $t$, the inequality
\[ 0\leq
\left(\displaystyle \frac{\lambda}{\sigma^2} \right) \int_0^{|u|}
r(x,t) f(x)dt\le \left(\displaystyle \frac{\lambda}{\sigma^2}
\right) |u|f\,, \]
and owing to the fact that $h(x,t)$ is an even function (because
$I_0 $ is even), it yields
\begin{equation}\label{ine}
h (x,u)+\left(\displaystyle \frac{\lambda}{\sigma^2} \right) |u|
f(x) \geq  \left(\displaystyle \frac{\lambda}{2\sigma^2} \right)
u^2 \,.
\end{equation}
Now, Young's inequality implies

{\color{black}
\[ 
\left(\frac{\lambda}{\sigma^2} \right) |u| f(x) \le \left(\frac{\lambda}{\sigma^2}
\right)  \displaystyle
\left[ \left( \frac{\epsilon}{2} \right) u^2+\left(\frac{1}{2\epsilon} \right)f(x)^2 \right]
\]
%
%\[ \left(\frac{\lambda}{\sigma^2} \right) |u| f(x) \le \epsilon\left(
%\frac{\lambda}{\sigma^2}
%\right)u^2+\frac1\epsilon\left(\frac{\lambda}{\sigma^2} \right)f(x)^2 \]
for any $\epsilon >0$. Thus, \eqref{ine}
becomes
\begin{equation}\label{ine0}
h (x,u)\ge \frac12 \left(
\frac{\lambda}{\sigma^2}
\right) \displaystyle \left[ \Big(1-\epsilon\Big)u^2-\frac1\epsilon f(x)^2\right],
\end{equation}
from where \eqref{linfty} follows choosing $\epsilon <1$}.

To go on, we need to know more features of the function
$\displaystyle s\mapsto\frac{r(s)}s$. Our starting point is
\eqref{amos3}. Indeed, letting $s\to0$ in \eqref{amos3}, it yields
$\displaystyle \lim_{s\to0}\frac{r(s)}s=\frac12$ and letting
$s\to+\infty$, we deduce $\displaystyle
\lim_{s\to+\infty}\frac{r(s)}s=0$. On the other hand,
\eqref{amos3} implies that the function $s\mapsto \frac{r(s)}s$ is
{\color{black}
(strictly)
}
decreasing in $[0,+\infty[$.
% Thus, we obtain
%\begin{equation}\label{amos4}
%\frac12>\frac{r(s)}s\downarrow0\,.
%\end{equation}

Next  let $w(x)$ be a positive critical point of
\[ h(x,t)=\displaystyle \left( \frac{\lambda }{2\sigma^2}\right) t^2
-\lambda \log I_0 \left( \frac{t f(x)}{\sigma^2} \right) \,, \]
then $h^\prime(x,w(x))=0$ and so $\displaystyle
w(x)=\displaystyle r\Big(\frac{f(x)w(x)}{\sigma^2}\Big)f(x)$. In other words,
\[ \frac{\sigma^2}{f(x)^2}=\displaystyle \frac{\displaystyle r\Big(\frac{f(x)w(x)}{\sigma^2}\Big)}{\displaystyle  \frac{f(x)w(x)}{\sigma^2}} \]%=\displaystyle \frac{r\big(f(x)w(x)/\sigma^2\big)}{f(x)w(x)/\sigma^2}\,.
According to \eqref{amos3}, it leads to the following dichotomy:
\begin{enumerate}
	\item If $0<f(x)^2\le 2\sigma^2$, then
	$\frac{\sigma^2}{f(x)^2}\ge\frac12$,  so that we cannot find a
	positive critical point. In this case $h(x,t)$ is convex and its
	minimum is attained at $0$.
	\item If $f(x)^2> 2\sigma^2$, then
  {\color{black}
	$0<\frac{\sigma^2}{f(x)^2}<\frac12$.
   As the function $s\mapsto \frac{r(s)}s$ is (strictly)
	decreasing, recall \eqref{amos3}, there exists a unique $s_f>0$
    satisfying
   \[
   \frac{\sigma^2}{f(x)^2}=\frac{r(s_f)}{s_f}\,.
   \]
   Choosing $w(x)$ such that $\frac{f(x)w(x)}{\sigma^2}=s_f$, we deduce that
   $w(x)>0$ and $h(x,t)$ has a critical point at
	$t=w(x)$.
   }
	Since $h(x,t)$ is negative in a neighbourhood of $0$ (as a consequence of $h''(x,0)<0$
	and $h' (x,0)=h(x,0)=0$) and
	$\lim_{t\to\pm\infty}h(x,t)=+\infty$, it follows that $h(x,t)$
	has, at least, a local minimum; wherewith that positive critical point  must be a local minimum. Therefore, $h(x,t)$ is an even function that, on $[0,+\infty[$, has the following
	profile: it is negative and decreasing in $[0,w(x)]$; it attains a global minimum at the point $w(x)$; from the point $w(x)$ on,  it is increasing; and
	goes to $+\infty$ as $t\to+\infty$.
\end{enumerate}
{\quad\rule{2mm}{2mm}\medskip}

\subsection{Coercitiveness and lower bound}\label{subsect_coerc}
In this Subsection we show that the energy minimization problem
related to the (formal) Euler--Lagrange equation in (\ref{hbes})
is coercive in $BV (\Omega )\cap L^2 (\Omega )$ because the energy
$H(u,f)$ defined in (\ref{hener}) is coercive in $L^2 (\Omega )$.
This shall be used to show that the energy $E_1 (u)$ has, at least, a
positive, non trivial minimum (provided that the datum is big enough).

Integrating  (\ref{ine0}) in $\Omega$, using definition
(\ref{hener}) and noticing that $r(x,0)=0$  we deduce:
\begin{multline*}
H (u,f)\geq \\ \Big(\frac12-\epsilon\Big)\left(\displaystyle
\frac{\lambda}{2\sigma^2} \right) \int_\Omega u^2 dx - \frac1\epsilon
\left(\displaystyle \frac{\lambda}{\sigma^2} \right) \int_\Omega
f^2 dx\,,
\end{multline*}
where $0< \epsilon < 1/2$,
and the functional $H(u,f)$ is coercive in $L^2 (\Omega )$. Then the energy functional
$E_p (u)$ in \eqref{ep} is coercive in $W^{1,p} (\Omega )\cap L^2
(\Omega )$ and $E_1 (u)$ (defined in \eqref{minFun}) is coercive in $BV (\Omega
)\cap L^2 (\Omega )$. These energies are also (uniformly) bounded
from below
\[ E_p (u)\geq - \frac1\epsilon \left(\displaystyle
\frac{\lambda}{\sigma^2} \right) \| f \|_2^2, \quad p\geq 1 \]
\subsection{Existence Result for the Approximating problems}\label{subsec:approx}
The analysis of problem \eqref{hbes} begins with the consideration of problems involving the $p-$Laplacian:
\begin{equation}\label{ap}
\left\{
\begin{array}{ll}
-\displaystyle \mbox{div} \left( |\nabla u |^{p-2} \nabla u\right)+h' (x, u)=0,&\mbox{in}\,\Omega , \\[0.5cm]
\displaystyle \left(  |\nabla u |^{p-2} \nabla u\right) \cdot n =0, & \mbox{on}\, \partial \Omega\,.
\end{array}
\right.
\end{equation}
Since we want let $p\to1$, it is enough to take $1<p< 2$. For such
$p$, the existence of a solution to \eqref{ap} is a standard
result although we have not found references for this specific
problem; so that we include its proof for the sake of
completeness. We are proving the following adaptation of Theorem
\ref{main}.

\begin{proposition}%\label{p-main}
	Let $1<p<2$ and  $\lambda >0$, $\sigma^2 \neq 0$ be given real parameters.
	
	For every non--negative $f\in L^2(\Omega)$, there exists a
	non--negative $u\in W^{1,p}(\Omega)\cap L^2(\Omega)$ which is a
	solution to problem \eqref{ap} and it is a {\em global} minimum of
	functional $E_p$.
\end{proposition}

{\sc Proof:} Consider the functional, written in terms of \eqref{hpr},
\begin{multline*}%\label{funcp}
E_p(u) = \\
\displaystyle \frac1p \int_\Omega |\nabla u |^{p}\, dx +\frac{\lambda}{2\sigma^2}\int_\Omega u^2 dx  -
\lambda\int_\Omega \log I_0\Big(\frac{uf}{\sigma^2}\Big) dx\,.
\end{multline*}
Since the Euler--Lagrange equation corresponding to the functional $E_p$ is \eqref{ap} and $E_p$ is differentiable, it is enough to
find a nonnegative  minimizer of $E_p$ in the space $W^{1,p}(\Omega)\cap L^2(\Omega)$.

The weakly lower--semicontinuity of $E_p$ can be obtained as follows. If $(u_n)_n$ is a sequence in  $W^{1,p}(\Omega)\cap L^2(\Omega)$ such that
\begin{gather*}
u_n\rightharpoonup u\,,\quad\hbox{weakly in }L^2(\Omega)\,;\\
\nabla u_n\rightharpoonup\nabla u\,,\quad\hbox{weakly in }L^p(\Omega;\mathbb{R}^N)\,;
\end{gather*}
then, due to the lower semicontinuity of the $p$--norm and the $2$--norm, it yields
\begin{gather*}
 \int_\Omega|\nabla u|^p\le \liminf_{n\to\infty}\int_\Omega|\nabla u_n|^p\\
\int_\Omega u^2\le \liminf_{n\to\infty}\int_\Omega u_n^2\,.
\end{gather*}
To pass to the limit in the remainder term, another consequence is in order, namely: the sequence $(fu_n)_n$ is weakly convergent in $L^1(\Omega)$, so that it is equi--integrable. Thus, it follows from the estimate
\[ \lambda\log I_0\Big(\frac{u_nf}{\sigma^2}\Big)\le \frac\lambda{\sigma^2}f|u_n| \]
that the sequence $\Big(\lambda\log I_0\Big(\frac{u_nf}{\sigma^2}\Big)\Big)_n$ is equi--integrable as well.
Moreover, applying the compact embedding of $W^{1,p}(\Omega)$ into $L^1(\Omega)$, we will assume that
\[u_n(x)\to u(x)\,, \quad\hbox{pointwise a.e. in }\Omega\,.  \]
This fact implies
\[\lambda\log I_0\Big(\frac{u_n(x)f(x)}{\sigma^2}\Big)\to \lambda\log I_0\Big(\frac{u(x)f(x)}{\sigma^2}\Big)\,, \ \hbox{a.e. in }\Omega\,. \]
By Vitali's Theorem we conclude that
\[ E_p(u)\le\liminf_{n\to\infty}E_p(u_n)\,. \]
On the other hand, we have already prove the coerciveness of $E_p$ previously in subsection \ref{subsect_coerc}. Therefore, there exists $u\in W^{1,p}(\Omega)\cap L^2(\Omega)$ which minimizes $E_p$.

Moreover, we may choose $u$ to be nonnegative. This feature is a consequence of being $h(x,s)$ an even function with respect to $s$, since
this fact induces
$
E_p(|u|)=E_p(u)
$
and so $|u|$ is a minimizer of $E_p$ as well. {\quad\rule{2mm}{2mm}\medskip}

\begin{remark}\rm
	Regarding uniqueness of problem \eqref{ap}, we point out that there always exists the trivial solution $u\equiv 0$. This solution may be unique if the datum is small enough (see section 4 below).
	
	Nevertheless, we are interested in uniqueness of positive solutions. When $p=2$, we may invoke  the results in \citep{Brezis1986} and, noting that the function $\displaystyle u\mapsto \frac{r(x,u)}u$ is decreasing, deduce that the positive solution to \eqref{ap} must be unique.
	Since $\displaystyle u\mapsto \frac{r(x,u)}{u^{p-1}}$ is not decreasing, this argument does not hold for $p<2$, so that we cannot presume that the positive solution we have found be unique.
\end{remark}

\section{Solving the model problem}\label{exis}
In this section we write rigorously the model equation formally
introduced in \eqref{hbes}. We shall prove the existence of a weak
(distributional) solution which is a global minimum of the energy
functional $E_1 (u)$ in \eqref{minFun}.

\subsection{Definition of solution for the model problem}\label{concept}
We shall say that $u\in BV (\Omega
)\cap L^2(\Omega)$ is a weak solution of problem \eqref{hbes} if
$h'(x,u)\in L^2(\Omega)$ and there exists a vector field
$\mathbf{z}\in L^{\infty} (\Omega , \mathbb{R}^N )$, with
$\|\mathbf{z} \|_\infty \leq 1$, such that
\begin{enumerate}
	\item $
	-\displaystyle  \mbox{div} (\mathbf{z}) + h' (x,u) =0 \qquad \mbox{in}\quad {\cal D}' (\Omega )
	$
	\item the equality $( \mathbf{z} , Du ) =|Du|$  holds in the sense of measures
	\item $[\mathbf{z}, n] =0$, ${\cal H}^{N-1}$--a.e. on the  boundary $\partial\Omega$.
\end{enumerate}

Roughly speaking, $\mathbf{z}$ plays the role of $\frac {Du}{|Du|}$.  The expressions
$( \mathbf{z} , Du )$ and $[\mathbf{z}, n]$ have sense thanks to the Anzellotti
theory (see \citep{Anzellotti1983} or \cite[Appendix C]{Andreu2004b}) which defines a Radon measure
$(\mathbf{z}, Dw)$, when $w\in BV (\Omega )\cap L^2(\Omega)$ and $\mbox{div} (\mathbf{z})\in L^2(\Omega)$, and provides the definition of a weakly trace on $\partial
\Omega$ to the normal component of $\mathbf{z}$, denoted by $[\mathbf{z}, n]$.
That Radon measure is defined, as a distribution, by the expression
\begin{equation}\label{defRM}
\langle (\mathbf z, Dw), \varphi \rangle=-\int_\Omega w\varphi\, \mbox{div} (\mathbf{z})\, dx-
\int_\Omega w\mathbf{z}\cdot \nabla \varphi\, dx\,,
\end{equation}
and its total variation satisfies the fundamental inequality
\begin{equation}\label{fund}
| (\mathbf z, Dw)|\le\|\mathbf z\|_\infty|Dw|\,.
\end{equation}
Furthermore, this theory also guarantees a Green's formula that
relates the function $[\mathbf{z}, n]$ and the measure
$(\mathbf{z}, Dw)$:
\begin{equation*}%\label{Green}
\int_{\Omega} w \ \mbox{div} (\mathbf{z}) \ dx + \int_{\Omega} (\mathbf{z}, Dw) =
\int_{\partial \Omega} [\mathbf{z}, n] w \ d {\mathcal H}^{N-1}\,.
\end{equation*}

Using this Green formula, we deduce a variational formulation of the solution to problem \eqref{hbes}, namely
\begin{equation}\label{sol}
\int_{\Omega} |Du| - \int_{\Omega} (\mathbf{z} , Dv) + \int_{\Omega} h'(x,u)(u-v) = 0\,,
\end{equation}
for all $v\in BV(\Omega)\cap L^2(\Omega)$.

This formulation allows us to show in which sense solutions to
problem \eqref{hbes} are critical points of the functional
$E_1=J_1+H$. In fact, it follows from \eqref{sol} that
\begin{eqnarray*}
 -\int_{\Omega} h'(x,u)(v-u) &=& \int_{\Omega} (\mathbf{z} ,
Dv)-\int_{\Omega} |Du|\\&\le& \int_{\Omega} |Dv|-\int_{\Omega} |Du|
\end{eqnarray*}
for all $v\in BV(\Omega)\cap L^2(\Omega)$. Hence,
\begin{equation*}%\label{multi}
-h^\prime(x,u)\in \partial J_1(u)\,.
\end{equation*}
\begin{remark}\label{rm31}\rm
	Observe that
	if we denote $\displaystyle
	F(u)=\lambda\log\Big(I_0\Big(\frac{fu}{\sigma^2}\Big)\Big)$, then
	we get that $F^\prime(u)$ lies in the subdifferential at $u$ of
	the convex functional defined by $\displaystyle v\mapsto
	\frac\lambda{2\sigma^2}\int_\Omega v^2+\int_\Omega|Dv|$.
\end{remark}

\subsection{A priori estimates}\label{ape}
We are proving that problem \eqref{hbes} has a solution $u$ for each $f\in L^2(\Omega)$. Moreover, we are getting $u\ge0$.

For $1<p< 2$, consider $u_p\in  W^{1,p}(\Omega)\cap L^2(\Omega)$ a nonnegative solution to the approximating problem
\begin{equation}\label{abes}
\left\{
\begin{array}{ll}
-\displaystyle \mbox{div} \left( |\nabla u_p |^{p-2} \nabla u_p\right)+h' (x, u_p)=0,&\mbox{in}\,\Omega , \\[0.5cm]
\displaystyle \left(  |\nabla u_p |^{p-2} \nabla u_p\right) \cdot n =0, & \mbox{on}\, \partial \Omega\,.
\end{array}
\right.
\end{equation}
The weak (variational) formulation of the boundary value problem \eqref{abes}, written in terms of \eqref{semi} and \eqref{ratio}, is:
\begin{eqnarray}\label{eq:1}
&&\displaystyle \frac{\lambda}{\sigma^2}\int_\Omega u_p vdx +  \displaystyle  \int_\Omega (|\nabla u_p |^{p-2} \nabla u_p )\cdot \nabla v dx \nonumber\\&=&
\displaystyle \frac{\lambda}{\sigma^2}\int_\Omega r(x,u_p)f vdx\,,
\end{eqnarray}
for all $ v\in W^{1,p} (\Omega )\cap L^2(\Omega)$.
Choosing $v=1$ we have the compatibility integral condition
\begin{equation}\label{comp}
\displaystyle \int_\Omega h' (x,u_p)dx =0
\end{equation}
i.e., $h' (x,u_p)$ has mean zero and we easily deduce a first estimate:
\[ \| u_p \|_1= \displaystyle \int_\Omega u_p dx =\displaystyle
\int_\Omega r(x,u_p)fdx \leq \displaystyle \int_\Omega fdx = \| f
\|_1 =M_1  \,. \]
\noindent We now use $v=u_p$ as a test function in the variational formulation
obtaining
\begin{eqnarray*}
&& \displaystyle \frac{\lambda}{\sigma^2}\int_\Omega u_p^2 dx +
 \displaystyle  \int_\Omega |\nabla u_p |^{p}   dx\\
 &=&
  \displaystyle
 \frac{\lambda}{\sigma^2}\int_\Omega r(x,u_p)f u_p dx\leq
 \displaystyle \frac{\lambda}{\sigma^2}\int_\Omega f u_p  \, dx\\
 &\leq& \displaystyle \frac{\lambda}{2\sigma^2}\displaystyle \left( \int_\Omega f^2 dx +\int_\Omega u_p^2 dx \right) \,,
\end{eqnarray*}
hence the uniform estimate
\begin{equation*}%\label{aga}
%\displaystyle \frac{1}{2\sigma^2} \| u_p \|_2^2 + \displaystyle \lambda \| \nabla u_p \|_p^p \leq M_2 = \displaystyle \frac{1}{2\sigma^2} \| f \|_2^2\,.
\displaystyle \lambda \| u_p \|_2^2 + \displaystyle 2\sigma^2  \|
\nabla u_p \|_p^p \leq  \displaystyle \lambda \| f \|_2^2 =M_2 \,.
\end{equation*}
It follows now from Young's inequality that
\[ \displaystyle \| u_p \|_2^2 +  \| \nabla u_p \|_1 \leq
\| u_p \|_2^2 + \displaystyle \frac1 p \| \nabla u_p \|_p^p+\frac{p-1}p|\Omega|\le \]
\[ \le \displaystyle \left(\frac{1}{\lambda}+ \frac{1}{2\sigma^2} \right)M_2 +|\Omega|=M_3\,. \]
Thus, $(u_p)_p$ is bounded in $BV(\Omega)\cap L^2(\Omega)$ and there exist $u\in BV(\Omega)\cap L^2(\Omega)$ and a subsequence, still denoted by $u_p$, satisfying
\begin{gather}
\label{conv1}    \nabla u_p \rightharpoonup Du\,,\quad\hbox{*--weakly as measures}\nonumber\\
\label{conv2}    u_p(x)\to u(x)\,,\quad\hbox{a.e. in }\Omega\nonumber\\
\label{conv3}   u_p\rightharpoonup u\,,\quad\hbox{weakly in }L^2(\Omega)\nonumber\\
\label{conv4}  u_p\to u\,,\quad\hbox{strongly in }L^r(\Omega)\quad\forall 1\le r<2
\end{gather}
We point out that $u\ge0$ due to being a pointwise limit of nonnegative functions.
We deduce from $u\in L^2(\Omega)$ that $h'(x,u)= (\lambda /\sigma^2)[u- r(x,u)f] \in L^2(\Omega)$ since $r(x,u)$ is bounded.
The boundedness of $(u_p)_p$ in $BV(\Omega)$ also implies that for every $q,\, 1\le q<p'$, we have

\begin{multline}\label{eq:2}
\displaystyle
\int_\Omega|\nabla u_p|^{(p-1) q}\, dx
\le \Big(\int_\Omega|\nabla u_p|^p\, dx\Big)^{(p-1) q/p} |\Omega|^{1-\frac{(p-1)q}{p}}\\
\\
\le M_3^{\frac{(p-1)q}{p}}
|\Omega|^{1-\frac{(p-1)q}{p}}\le M_3+|\Omega|\,.
\end{multline}
So, for any $q>1$ fixed, the sequence $|\nabla
u_p|^{p-2}
\nabla u_p$ is bounded in
$L^q(\Omega;\mathbb R^N)$ and then there exists $\mathbf z_q\in L^q(\Omega;\mathbb R^N)$ such
that, up to subsequences,
$$
|\nabla u_p|^{p-2}\nabla u_p\rightharpoonup \mathbf z_q\quad\hbox{in }
L^q(\Omega;\mathbb R^N)\quad\hbox{for all  } 1\le q<+\infty\,. $$
Moreover, by a diagonal argument we can find a limit $\mathbf z$ that does not depend on
$q$, that is
\begin{equation}\label{eq:3}
|\nabla u_p|^{p-2}\nabla u_p \rightharpoonup \mathbf z\quad\hbox{in }
L^q(\Omega;\mathbb R^N)\ \hbox{ for  } 1\le q<+\infty\,.
\end{equation}
Now by \eqref{eq:2} we deduce
\[\||\nabla u_p|^{p-2} \nabla u_p\|_{L^q(\Omega;\mathbb R^N)}\le (M_3+|\Omega|)^{1/q} \]
for  $1\le q<+\infty$ and for $
p\in]1,q^\prime[\,$. Therefore, by lower semicontinuity of the
norm, we have
\begin{equation*}
\|\mathbf z\|_{L^q(\Omega;\mathbb R^N)}\le (M_3+|\Omega|)^{1/q}\quad\hbox{for
	all  } 1\le q<+\infty\,.
\end{equation*}
Letting $q\to \infty$, we get that
$\mathbf z\in L^\infty(\Omega;\mathbb R^N)$ and
\begin{equation*}
%\label{eq:4}
\|\mathbf z\|_{L^\infty(\Omega;\mathbb R^N)}\le1\,.
\end{equation*}

\subsection{Checking that function $u$ is a solution to the model problem \eqref{hbes}}

We have to see that $u$ satisfies the requirements of our
definition (see Subsection \ref{concept} above).

Taking $v=\varphi\in C_0^\infty(\Omega)$ in \eqref{eq:1} and
letting $p\to1$, it yields
\[ \displaystyle \frac{\lambda}{\sigma^2}\int_\Omega u \varphi dx +
\displaystyle   \int_\Omega \mathbf z\cdot \nabla \varphi dx =
\displaystyle \frac{\lambda}{\sigma^2}\int_\Omega r(x,u)f \varphi
dx\,, \]
so that our equation holds in the sense of distributions.

Once we have proved 1 in the definition of solution, we proceed to see 2 and 3. To begin with 2, consider $\varphi\in C_0^\infty(\Omega)$ such that $\varphi\ge0$. Taking $u_p\varphi$ as test function in \eqref{eq:1}, we obtain
\begin{eqnarray}\label{eq:5}
&&\frac{\lambda}{\sigma^2}\int_\Omega u_p^2\varphi dx+\int_\Omega \varphi|\nabla u_p|^pdx\\&&+
\int_\Omega u_p |\nabla u_p|^{p-2}\nabla u_p\cdot \nabla \varphi dx=\frac{\lambda}{\sigma^2}\int_\Omega r(x,u_p)fu_p\varphi dx\,.
\nonumber
\end{eqnarray}
We are studying each term in \eqref{eq:5} to let $p\to 1$. We apply Fatou's Lemma in the first term. In the second, we use Young's inequality and the lower semicontinuity of the total variation as follows:
\begin{eqnarray*}
\int_\Omega \varphi|\nabla u_p|&\le&\liminf_{p\to1}
\int_\Omega \varphi|\nabla u_p|\,dx\\
&\le&\liminf_{p\to1}\Big(\frac1p\int_\Omega \varphi|\nabla u_p|^pdx+\frac{p-1}p\int_\Omega\varphi\,dx\Big)\\
&=&\liminf_{p\to1}\int_\Omega \varphi|\nabla u_p|^pdx\,.
\end{eqnarray*}
Third term is handled using \eqref{conv4} and \eqref{eq:3}. In the right hand side is enough to have in mind that $r$ is bounded. Therefore, \eqref{eq:5} becomes
\begin{eqnarray*}
&&\frac{\lambda}{\sigma^2}\int_\Omega u^2\varphi dx+\int_\Omega \varphi|D u|+
\int_\Omega u \mathbf z\cdot \nabla \varphi dx\\&&\le \frac{\lambda}{\sigma^2}\int_\Omega r(x,u)fu\varphi dx\,.
\end{eqnarray*}
Taking into account that our equation holds in the sense of distributions and simplifying, we may write this inequality as
\begin{equation*}
\int_\Omega \varphi|D u|+
\int_\Omega u \mathbf z\cdot \nabla \varphi \,dx\le-\int_\Omega u\varphi\, \mbox{div}\mathbf z\, dx\,.
\end{equation*}
By \eqref{defRM}, this is just
\[ \int_\Omega \varphi|D u|\le \langle (\mathbf z, Du),\varphi\rangle\,, \]
that is, $|Du|\le(\mathbf z,Du)$ as measures. The reverse inequality is a consequence of \eqref{fund}. Hence, 2 is seen.

It only remains to prove 3. To this end, consider $v\in W^{1,2}(\Omega)$ in \eqref{eq:1} and
take limits as $p$ goes to $1$. It yields
\[ \frac{\lambda}{\sigma^2}\int_\Omega uv\,dx+ \int_\Omega \mathbf z\cdot\nabla v\, dx=\frac{\lambda}{\sigma^2}\int_\Omega r(x,u)fv\, dx\,. \]
Using the equality $\displaystyle\frac{\lambda}{\sigma^2}(u- r(x,u)f)=\,\mbox{div\,}\mathbf z$, it follows that
\[ \int_\Omega v\,\mbox{div\,}\mathbf z\,dx+ \int_\Omega \mathbf z\cdot\nabla v\, dx=0\,, \]
so that Green's formula implies
\[ \int_{\partial\Omega}v[\mathbf z, n]\, d \mathcal H^{N-1}=0\,. \]
By a density argument, this leads to $[\mathbf z, n]=0$ $\mathcal H^{N-1}$--a.e. on
$\partial\Omega$.  {\quad\rule{2mm}{2mm}\medskip}

\begin{remark}\label{compc}\rm
	
	We explicitly point out that the compatibility condition
	\eqref{comp} also holds for the solution $u$ to problem
	\eqref{hbes}. To check this fact, it is enough to multiply
\[ 	-\hbox{div\,}(\mathbf z)+h'(x,u)=0 \]
	by a constant function and apply Green's formula. Then we get
	\[ \int_\Omega h'(x,u) dx=0\,. \]
	The same condition can be deduced letting $p\to1$ in \eqref{comp}
	since $|h'(x,u_p)|\le C(|u_p|+f)$ and $u_p\to u$ strongly in
	$L^1(\Omega)$.
\end{remark}

\subsection{Function $u$ is a global minimizer of functional $E_1$}

We will prove that the nonnegative
{\color{black}
function $u$ considered in Subsection 3.2, which we have shown is a solution to problem \eqref{hbes}
in Subsection 3.3,
}
satisfies
\begin{equation*}%\label{glmin}
E_1(u)\le E_1(v)\,,\qquad\hbox{for all }v\in BV(\Omega)\cap L^2(\Omega)\,.
\end{equation*}
To see it, we use several stages.

{\bf Step 1.-} To begin with, assume that $v\in W^{1,2}(\Omega)$.
Observe first that the interpolation inequality implies
\[ \frac1p\int_\Omega|\nabla v|^p\, dx\le\frac1p \|\nabla
v\|_2^{2(p-1)}\|\nabla v\|_1^{p-2(p-1)} \]
for all $1<p<2$. Thus,
\[ \limsup_{p\to1}\frac1p\int_\Omega|\nabla v|^p\, dx\le
\int_\Omega|\nabla v|\, dx\,. \]
On the other hand, as a consequence of Young's inequality, we have
\[ \int_\Omega|\nabla v|\, dx\le \frac1p \int_\Omega|\nabla v|^p\,
dx+\frac{p-1}p|\Omega|\,, \]
for all $1<p<2$; so that
\[ \int_\Omega|\nabla v|\, dx\le \liminf_{p\to1}\frac1p
\int_\Omega|\nabla v|^p\, dx\,. \]
Hence, the conclusion is
\[ \int_\Omega|\nabla v|\, dx= \lim_{p\to1}\frac1p
\int_\Omega|\nabla v|^p\, dx\,, \]
that is
\begin{equation}\label{limp}
E_1(v)= \lim_{p\to1}E_p(v)\,.
\end{equation}

Since $u_p$ is a minimizer of $E_p$ and $v\in W^{1,p}(\Omega)\cap
L^2(\Omega)$, we obtain
\[ E_p(u_p)\le \frac1p\int_\Omega|\nabla v|^p\, dx+\int_\Omega
h(x,v)\, dx\,, \]
for all $1<p<2$. On account of \eqref{limp}, using the
lower--semicontinuity of functional $E_1$ and Young's inequality
we deduce that
\begin{eqnarray*}
E_1(u)&\le& \liminf_{p\to1}E_1(u_p)\le
\liminf_{p\to1}\Big(E_p(u_p)+\frac{p-1}p|\Omega|\Big)
\\&\le&
\lim_{p\to1}E_p(v)=E_1(v)\,.
\end{eqnarray*}

{\bf Step 2.-} Assume now that $v\in W^{1,1}(\Omega)\cap
L^2(\Omega)$ satisfies $v\big|_{\partial\Omega}\in
W^{1/2,2}(\partial\Omega)$. Then there exists $w\in
W^{1,2}(\Omega)$ such that
$v\big|_{\partial\Omega}=w\big|_{\partial\Omega}$ and so $v-w\in
W_0^{1,1}(\Omega)\cap L^2(\Omega)$. Thus, there exists a sequence
$(v_n)_n$ in $C_0^\infty(\Omega)$ such that
\begin{gather*}
v_n+w\to v\,,\quad\hbox{strongly in }W^{1,1}(\Omega)\,;\\
v_n+w\to v\,,\quad\hbox{strongly in }L^{2}(\Omega)\,.
\end{gather*}
Since Step 1 provides us
\[ E_1(u)\le E_1(v_n+w)\,,\hbox{for all }n\in \mathbb N\,, \]
it follows that
\[ E_1(u)\le E_1(v)\,. \]
{\bf Step 3.-} Consider the general case: $v\in BV(\Omega)\cap
L^2(\Omega)$. Some approximation sequences of $v$ are in order.
First (see \cite[Theorem 3.9]{Ambrosio2000} and \cite[Remark 2.12]{Giusti1984}) there exists a sequence $(v_n)_n$ in $C^\infty(\Omega)\cap
W^{1,1}(\Omega)\cap L^2(\Omega)$ such that
\begin{gather*}
v_n\to v\,,\quad\hbox{strongly in }L^{2}(\Omega)\,,\\
\int_\Omega|\nabla v_n|\to \int_\Omega|Dv|\,,\\
v_n\big|_{\partial\Omega}=v\big|_{\partial\Omega}\,,\quad\hbox{for
	all }n\in \mathbb N\,.
\end{gather*}
On the other hand, given $v\big|_{\partial\Omega}\in
L^1(\partial\Omega)$, we may find a sequence $(\varphi_n)_n$ in
$W^{1/2,2}(\partial\Omega)$ satisfying
\begin{equation*}
\varphi_n\to v\big|_{\partial\Omega}\,,\quad\hbox{strongly in }L^{1}(\partial\Omega)\,.
\end{equation*}

For each $n\in\mathbb
N$, we apply \cite[Lemma 5.5]{Anzellotti1983} to get $w_n\in C(\Omega)\cap
W^{1,1}(\Omega)\cap L^2(\Omega)$ such that
\begin{gather*}
\int_\Omega|\nabla w_n|\, dx<\int_{\partial\Omega}|\varphi_n-v|\, d\mathcal H^{N-1}+\frac1n\,,\\
\int_\Omega|w_n|^2\, dx<\frac1n\,,\\
w_n\big|_{\partial\Omega}=\varphi_n-v\big|_{\partial\Omega}\,,\quad\hbox{for
	all }n\in \mathbb N\,.
\end{gather*}

Summing up, we have
\begin{enumerate}
	\item $w_n+v_n\in C(\Omega)\cap W^{1,1}(\Omega)\cap
	L^2(\Omega)$ for  all $n\in \mathbb N$;
	\item $w_n+v_n\to v$, strongly in $L^2(\Omega)$;
	\item $(w_n+v_n)\big|_{\partial\Omega}=\varphi_n\in W^{1/2,2}(\partial\Omega)$ for  all $n\in \mathbb
	N$.
\end{enumerate}
Moreover, since
\begin{multline*}
\int_\Omega|\nabla (w_n+v_n)|\, dx\le\int_\Omega|\nabla w_n|\,
dx+\int_\Omega|\nabla v_n|\, dx\\
<\int_{\partial\Omega}|\varphi_n-v|\, d\mathcal
H^{N-1}+\frac1n+\int_\Omega|\nabla v_n|\, dx
\end{multline*}
and $\varphi_n\to v\big|_{\partial\Omega}$ strongly in
$L^{1}(\partial\Omega)$, it follows that
\[ \limsup_{n\to\infty} \int_\Omega|\nabla (w_n+v_n)|\, dx\le
\int_\Omega|Dv|\,. \]
The lower semicontinuity of the total variation now leads to
\[ \lim_{n\to\infty} \int_\Omega|\nabla (w_n+v_n)|\, dx=
\int_\Omega|Dv|\,. \]
Therefore,
\[ E_1(w_n+v_n)\to E_1(v)\,. \]
Finally, by Step 2, we already get
\[ E_1(u)\le E_1(w_n+v_n)\,,\quad\hbox{for all }n\in \mathbb N\,. \]
Letting $n$ go to $\infty$, we see that $E_1(u)\le E_1(v)$. Since
this fact holds for all $v\in BV(\Omega)\cap L^2(\Omega)$, we are done.
\section{Remarks and Properties of the Problem}

\subsection{Summability of the solutions}
We are interested in dealing with bounded data $f$. In this case,
the solution we find is also bounded. More generaly, we will see
in this remark that if $f\in L^q(\Omega)$, with $q>N$, then the
solution $u$ is bounded. It is enough to check that an
$L^\infty$--estimate holds on the approximate solutions $u_p$.
Since $q > N$, then $\frac{N}{q' (N-1)} > 1$. Fix $p_0$, such that
$1 < p_0 < \frac{N}{q' (N-1)}$, and take $p$ such that $1 < p \leq
p_0$. For any $k >0$ consider the real function
{\color{black}
$G_k(s) := (s - k)^+$, $s \ge 0$.
}
Taking $G_k(u_p)$ as test function in \eqref{eq:1}, we get
\begin{eqnarray*}
&&\frac\lambda{\sigma^2}\int_{\Omega} u_p G_k(u_p) \, dx+ \int_{\Omega} \vert \nabla G_k(u_p) \vert^p\, dx \\&\leq&  \frac\lambda{\sigma^2}\int_{\Omega} f r(x,u_p)  G_k(u_p)\, dx\,.
\end{eqnarray*}
Disregarding a nonnegative term and applying $r(x,u_p)\le1$,  H\"older's inequality leads to
\begin{eqnarray*}
\int_{\Omega} \vert \nabla G_k(u_p) \vert^p\, dx &\leq& \frac\lambda{\sigma^2}\int_{\Omega} f  G_k(u_p)\, dx
\\&\le& \|f\|_q\Big(\int_\Omega|G_k(u_p)|^{q^\prime}\, dx\Big)^{1/q^\prime}\,.
\end{eqnarray*}
This is the starting point for using the Stampacchia technique and get an $L^\infty$--estimate. Just be careful to check that the various constants appearing in the calculations do not depend on $p$.
Details can be found at \cite[Theorem 3.5, Step 3]{Mazon2013}.

\bigskip

Furthermore, if $f \in L^\infty(\Omega)$, we may clarify a little more the situation by seeing the estimate $\|u\|_\infty\le \|f\|_\infty$.
This inequality makes explicit and extends the statement 8 of
\cite[Theorem 1]{Getreuer2011}.

Taking $u_p^q$, with $q>1$ large enough, as test function and dropping a nonnegative term, we obtain
\[ \frac\lambda{\sigma^2}\int_{\Omega} u_p^{q+1} \, dx \le  \frac\lambda{\sigma^2}\int_{\Omega} f r(x,u_p)  u_p^q\, dx
\le  \frac\lambda{\sigma^2}\int_{\Omega} f   u_p^q\, dx\,. \]
It follows from H\"older's inequality that
\begin{eqnarray*}
\int_{\Omega} u_p^{q+1} \, dx
&\le&  \int_{\Omega} f   u_p^q\, dx\\
&\le& \Big(\int_{\Omega} f^{q+1} \, dx \Big)^{1/(q+1)}\Big(\int_{\Omega} u_p^{q+1} \, dx \Big)^{q/(q+1)}
\,,
\end{eqnarray*}
and so
\[ \Big(\int_{\Omega} u_p^{q+1} \, dx  \Big)^{1/(q+1)}
\le \Big(\int_{\Omega} f^{q+1} \, dx \Big)^{1/(q+1)}\,. \]
Letting $q\to\infty$, it yields $\|u_p\|_\infty\le\|f\|_\infty$
for all $1<p<2$, and recalling that $u$ is the pointwise limit of
$u_p$, we are done.

\subsection{Uniqueness}\label{uniqu} We will prove that if the function $t\mapsto
h^\prime(x,t)$ is increasing, then there exists at most a solution
to \eqref{hbes}.

\begin{figure}[h]\label{fig4}
	\begin{center}
		{\includegraphics[width=0.9\columnwidth]{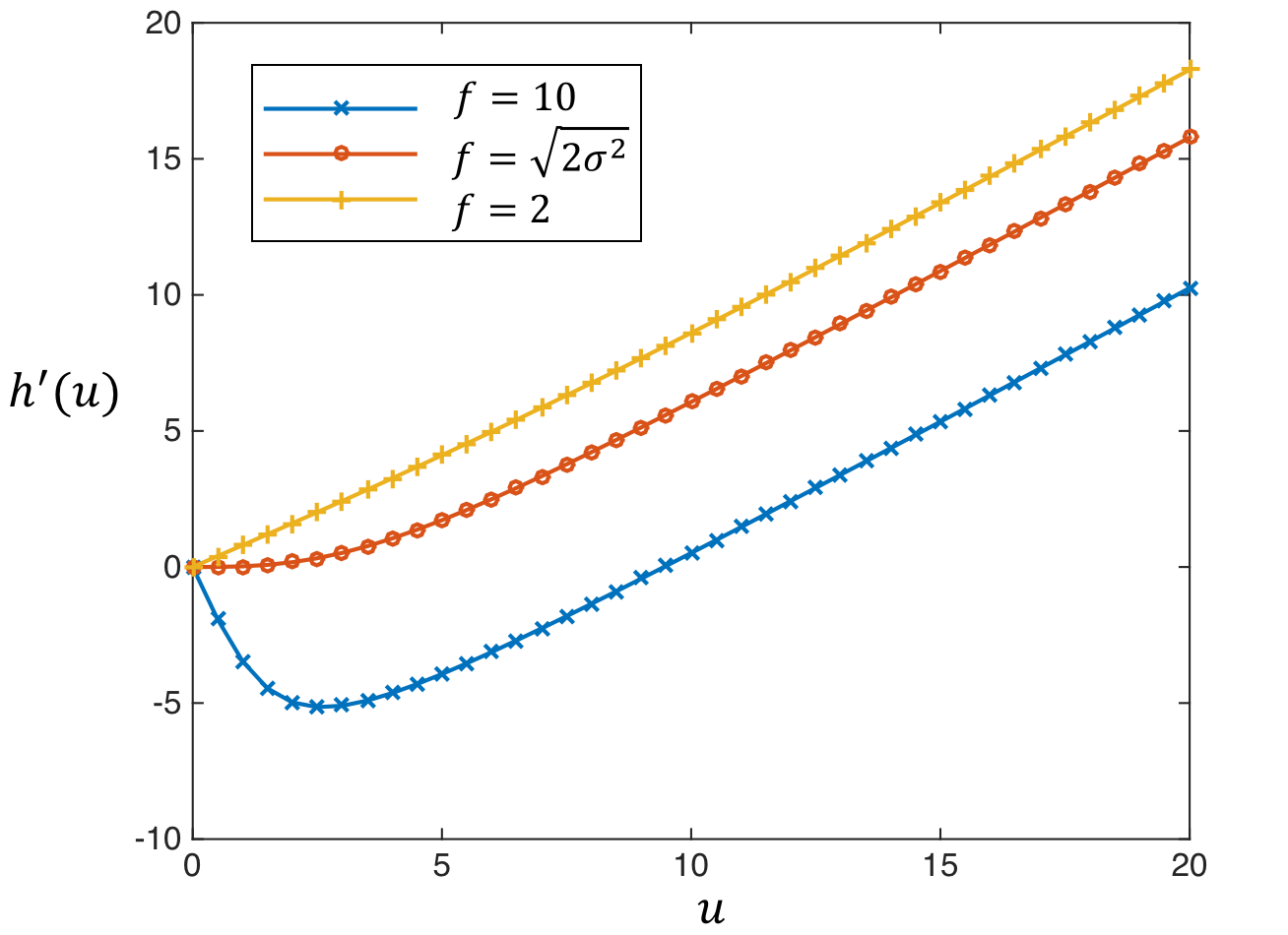}}
	\end{center}
	\caption{Profile of $h^{'}(x,u)$ for fixed $x\in\Omega$ and the parametric values $\lambda =10$, $\sigma^2 =10$  for different, constant  values of the data: $f=2$, $f=f^* =\sqrt{20}$ and $f=10$.
		A limit behavior is obtained when $f=f^* =\sqrt{20} =\sqrt{ 2\sigma^2}$. For $f\leq f^* $ we have uniqueness of the trivial solution.
		For $f>f^* $  we have $f^2 >2\sigma^2$ and the corresponding profile is negative in a neighborhood
		of $s=0$. Notice that when $u$ is small, $h' <0$ and $h'$ behaves as a reactive term (a source) in the Euler-Lagrange equation. When $u$ is sufficiently big $h' >0$ and $h'$ define an absorption term (a sink) in the equation.}
\end{figure}
{\sc Proof:}
Assume, to get a contradiction, that $u_1$ and $u_2$ are two
solutions to \eqref{hbes} in the sense of the definition stated in
Subsection \ref{concept} above. Denote by $\mathbf z_1$ and
$\mathbf z_2$ the respective vector fields. It follows that
\begin{equation*}
-\hbox{div\,}\mathbf z_i+h'(x,u_i)=0\,,\qquad i=1,2\,;
\end{equation*}
in the sense of distributions. Multiply both equations by
$u_1-u_2$, use Green's formula, recall the second condition in the
Definition of solution to problem \eqref{hbes} and substract one
expresion from the other to obtain
\begin{eqnarray*}
&&\int_\Omega |Du_1|-(\mathbf z_2, Du_1)+\int_\Omega |Du_2|-(\mathbf
z_1, Du_2)\\
&+&\int_\Omega(h'(x,u_1)-h'(x,u_2))(u_1-u_2)\, dx=0\,.
\end{eqnarray*}
The three terms are nonnegative since $(\mathbf z_i , Du_j)\le
\|z_i\|_\infty |Du_j|\le  |Du_j| $, for $i,j=1,2$, and the
function $t\mapsto h^\prime(x,t)$ is increasing. Hence, they must
vanish; in particular,
\[ \int_\Omega(h'(x,u_1)-h'(x,u_2))(u_1-u_2)\, dx=0 \]
and $h'$ increasing implies $u_1\equiv u_2$, as desired. {\quad\rule{2mm}{2mm}\medskip}

\subsection{Non trivial solutions} We have already commented that there
always exists a trivial solution $u\equiv 0$. On the other hand,
$0\le f\le \sqrt{2\sigma^2}$ implies that $h'(x,s)$ is increasing
with respect to $s$ and, as a consequence of the uniqueness result
of the previous subsection, there is no other solution aside from the
trivial one. Nevertheless, we are interested in the case when
$f\in L^\infty(\Omega)$ is a.e. above this threshold  and in finding non trivial
solutions.

{\color{black}
Although constant data are unrealistic, we study them to get nontrivial solutions.
In this Subsection, }
we are showing that if the datum is constant $f(x)=\mu$ and
$\mu>\sqrt{2\sigma^2}$, then the solution is constant and
non trivial.
{\color{black}
It is worth remarking that we obtain uniqueness of positive solutions for constant data.
}
{\color{black}
In the next Subsection, we will derive a criterion on the datum to obtain nontrivial solutions.
}

Considering \eqref{hpr},  we define the function
\[ \Gamma (\mu ,t)=\displaystyle \left( \frac{\lambda }{2\sigma^2}\right) t^2 -\lambda \log I_0 \left( \frac{t \mu}{\sigma^2}
\right) \]
{\color{black} which is related to the function $h(x,u)$ setting $h(x,u)=h(f(x),u(x))$. Fixed $x\in \Omega $ we have
$\Gamma (\mu , t)= h(f(x),u(x))$ and we can use the results in the proof Lemma \eqref{lemma2}, condition 2}.  Computing its derivative we have, $\forall\, \mu >0,\, t>0$,
\begin{equation}\label{decr}
\Gamma_\mu (\mu ,t)=-  \frac{\lambda t}{\sigma^2} \displaystyle
\left[  \frac{I_1 (t \mu /\sigma^2 )}{I_0 (t \mu /\sigma^2 )}
\right] <0
\end{equation}
and $\Gamma (\mu ,t)$ is decreasing with respect to $\mu$. Owing
to $\mu>\sqrt{2\sigma^2}$,  the function $\Gamma (\mu ,t)$ attains
a negative minimum at a positive point, say $t=\gamma$ {\color{black} (see the end of the proof in Lemma \eqref{lemma2} , condition 2)}. Then,
fixed $\mu$, $\gamma = \argmin \, \Gamma (\mu ,t)$ satisfies
$\Gamma_t (\mu ,t) =0$ which is
\begin{equation}\label{uniq}
\displaystyle
\gamma=\left[\frac{I_1(\gamma\mu/\sigma^2)}{I_0(\gamma\mu/\sigma^2)}\right]\mu\,.
\end{equation}
Actually, there is just a positive point $\gamma$ satisfying
\eqref{uniq}; to see this it is enough to check that
$s_\mu=\frac{\mu\gamma}{\sigma^2}$ is the unique solution to
problem
\[ s_\mu=\frac{I_1(s_\mu)}{I_0(s_\mu)}\frac{\mu^2}{\sigma^2}=r(s_\mu)\frac{\mu^2}{\sigma^2} \]
and this fact is a consequence of being the function $s\mapsto
\frac{r(s)}s$ decreasing in $[0,+\infty[$ (see \eqref{amos3}).
Thus $\gamma$ is given by \eqref{uniq} and satisfies
\begin{equation}\label{uniq2}
\Gamma (\mu ,\gamma)<\Gamma (\mu ,0)=0\quad \hbox{and} \quad0< \gamma <\mu\,.
\end{equation}

{\color{black}
Taking $u(x)=\gamma$ for all $x\in\Omega$, it yields that $u$ is
the unique minimizer of the functional $E_1$. Indeed, if $v\in
BV(\Omega)\cap L^2(\Omega)$, then $h(x,u)\le h(x,v)$ and
$h(x,u)= h(x,v)$ only when $v(x)=\gamma$ a.e., so that
\[ E_1(u)=\int_\Omega h(x,u)\, dx\le \int_\Omega |Dv|+\int_\Omega
h(x,v)\, dx=E_1(v) \]
and $E_1(u)=E_1(v)$ only when $u=v$.
}

\subsection{Comparing with constant functions}
Using the same notation of the above subsection, we may go further and prove that
{\color{black}
$0\le \mu\le f(x)$ implies $\gamma\le u(x)$
a.e.in $\Omega$, where $\gamma\ge 0$ minimizes
$\Gamma(\mu,t)$. We also assume that
$\mu>\sqrt{2\sigma^2}$, otherwise $\gamma=0$ and the inequality becomes obvious.
}

We begin by claiming that, for almost all $x\in\Omega$,
\begin{equation}\label{claim0}
\hbox{function  }t\mapsto h(x,t) \hbox{ is  (strictly) decreasing in  } [0,\gamma]
\end{equation}
{\color{black}
and
\begin{equation}\label{claim0.0}
\hbox{function  }t\mapsto h(x,t) \hbox{ is  (strictly) increasing in  } [\gamma_2,+\infty[\,.
\end{equation}
In both cases,
}
we will use that functions
\begin{equation}\label{claim1}
s\mapsto \frac{s I_0(s)}{I_1(s)}\quad\hbox{and}\quad s\mapsto
\frac{I_1(s)}{I_0(s)}\quad\hbox{ are increasing}
\end{equation}
and these facts are derived from \eqref{amos3}. Notice that, for
almost all $x\in\Omega$, the positive minimum $w(x)$ of $h(x,t)$
satisfies
\begin{equation}\label{claim2}
w(x)=\left[\frac{I_1\left(f(x)w(x)/\sigma^2\right)}{I_0\left(f(x)w(x)/\sigma^2\right)}\right]f(x)\,.
\end{equation}
It follows that
\[ s_f(x)=\left[\frac{I_1(s_f(x))}{I_0(s_f(x))}\right]\frac{f(x)^2}{\sigma^2}\,, \]
where $s_f(x)=\frac{f(x)w(x)}{\sigma^2}$. As seen in the previous
subsection, a similar identity holds for the positive minimum
{\color{black}
$\gamma$ of $\Gamma(\mu,t)$:
\[ s_{\mu}=\left[\frac{I_1(s_{\mu})}{I_0(s_{\mu})}\right]\frac{\mu^2}{\sigma^2}\,, \]
where $s_{\mu}=\frac{\mu\gamma}{\sigma^2}$.
Hence, by \eqref{claim1}, $\mu\le f(x)$ implies $s_{\mu}\le
s_f(x)$ a.e. and so $\mu\gamma\le f(x)w(x)$ a.e.
Going back to
\eqref{claim2}, for almost all $x\in\Omega$, we have
\begin{eqnarray*}
w(x)^2&=&\left[\frac{I_1\left(f(x)w(x)/\sigma^2\right)}{I_0\left(f(x)w(x)/\sigma^2\right)}\right]f(x)w(x)
\\&\ge&
\left[\frac{I_1\left(f(x)w(x)/\sigma^2\right)}{I_0\left(f(x)w(x)/\sigma^2\right)}\right]\mu\gamma\\
&\ge&
\left[\frac{I_1\left(\mu\gamma/\sigma^2\right)}{I_0\left(\mu\gamma/\sigma^2\right)}\right]\mu\gamma=\gamma^2\,,
\end{eqnarray*}
where the last inequality is due to \eqref{claim1}. Therefore, we
have seen that $w(x)\ge\gamma$ a.e. Finally, since $h(x,\cdot)$
is decreasing in $[0,w(x)]$ for almost all $x\in \Omega$, it
yields that $h(x,\cdot)$ is decreasing in $[0,\gamma]$ for
almost all $x\in \Omega$ and \eqref{claim0} is proved.
}
{\color{black}
The second
claim follows using a similar argument.
}

{\color{black}
Now we turn to check that $u(x)\ge\gamma$ a.e. Since $u$ is a global
minimizer of functional $E_1$, it follows that
\begin{eqnarray*}
&&\int_\Omega|Du|+\int_\Omega h(x,u)\, dx\\
&\le&\int_\Omega|D(u+(\gamma-u)^+) |+\int_\Omega h(x,u+(\gamma-u)^+)\, dx\\
&\le&\int_{\Omega}|Du|+\int_{\{u\ge\gamma\}} h(x,u)\, dx +\int_{\{u<\gamma\}} h(x,\gamma)\, dx\,.
\end{eqnarray*}
Simplifying and dropping the nonnegative gradient term, we obtain
\[ \int_{\{u<\gamma\}}h(x,u)\, dx\le
\int_{\{u<\gamma\}}h(x,\gamma)\, dx\,. \]
Applying now our first claim \eqref{claim0}, we deduce that
$h(x,\gamma)< h(x,u)$ a.e. in $\{u<\gamma\}$. Therefore,
$|\{u<\gamma\}|=0$, that is $u(x)\ge\gamma$ a.e. in $\Omega$.
}

{\color{black}
Starting from the inequality
\begin{eqnarray*}
&&\int_\Omega|Du|+\int_\Omega h(x,u)\, dx\\
&\le& \int_\Omega|D(u-(u-\gamma_2)^+)|+\int_\Omega h(x,u-(u-\gamma_2)^+)\,
dx\,,
\end{eqnarray*}
it follows that
\[ \int_{\{u>\gamma_2\}}h(x,u)\, dx\le
\int_{\{u>\gamma_2\}}h(x,\gamma_2)\, dx\,. \]
and our second claim \eqref{claim0.0} implies that
$u(x)\le\gamma_2$ a.e. in $\Omega$.
}
\subsection{The minimum is decreasing with respect to the datum} In
this remark we will make explicit the dependence on the data. To
this end, we stand our functional for $E^f_1$.

Let $f_i\in L^2(\Omega)$, $i=1,2$, be two data and denote by $u_i$
the corresponding function where the minimum of $E^{f_i}_1$ is
attained. We will show that $f_1\le f_2$ implies
$E_1^{f_1}(u_1)\ge E_1^{f_2}(u_2)$.

Since $f_1(x)\le f_2(x)$ implies
$H(v,f_1)\ge H(v,f_2)$ for all $v\in BV(\Omega)\cap L^2(\Omega)$,
recall \eqref{decr}, it follows that
\begin{eqnarray*}
E_1^{f_2}(u_2)&\le& E_1^{f_2}(u_1)=\int_\Omega|Du_1|+H(u_1,f_2)\\
&\le& \int_\Omega|Du_1|+H(u_1,f_1)=E_1^{f_1}(u_1)\,.
\end{eqnarray*}

Combining this fact with the previous subsection and having in mind
\eqref{uniq2}, we get that $f(x)\ge\mu>\sqrt{2\sigma^2}$ implies
\[ E_1^f(u)\le E_1^\mu(\gamma)=\Gamma(\mu,\gamma)|\Omega|<0\,. \]

\subsection{Resolvents of the subdifferential}\label{ressub}
With a view to the numerical resolution of problem \eqref{hbes}, we now consider some properties
of the resolvents of the sub--differential of
a (possibly) quadratically  perturbed Total Variation energy functional.

It is well--known that subdifferentials of convex functions have nonexpansive resolvents.
Thanks to the characterization of the subdifferential of the Total Variation appearing in \citep{Andreu2004b},
we may make explicit this feature in our case. Indeed, fix $\alpha\ge0$ and set
\begin{multline}\label{Guno}
G_1 (u)=
\left\{
\begin{array}{ll}
\displaystyle
\int_\Omega  |D u|  + \alpha \| u\|_2^2\,,&\hbox{if }u\in BV(\Omega)\cap L^2(\Omega)\,;\\[5mm]
+\infty\,,&\hbox{if }u\in L^2(\Omega)\backslash BV(\Omega)\,.
\end{array}
\right.
\end{multline}
Using \cite[Lemma 2.4]{Andreu2004b}, it yields that $u\in (I+c\,\partial G_1)^{-1}(f)$, with $c>0$, if and only if $u$ is a solution to
\begin{equation}\label{resolv}
\left\{
\begin{array}{ll}
\displaystyle u+c\alpha u-c\,\hbox{div\,}\Big(\frac{Du}{|Du|}\Big)=f\,,&\hbox{in }\Omega\,;\\[5mm]
\displaystyle \Big(\frac{Du}{|Du|}\Big)\cdot n=0\,,&\hbox{on }\partial\Omega\,.
\end{array}
\right.
\end{equation}
We point out that this problem has a unique solution (just follow the arguments in subsection \ref{uniqu}).

Consider now $u_i$ solution to problem \eqref{resolv} with datum $f_i$, $i=1,2$. In other words, we have
$u_i= (I+c\,\partial G_1)^{-1}(f_i)$, $i=1,2$.
Then there exist $\mathbf z_i\in L^\infty(\Omega;\mathbb R^N)$ satisfying the requirements of Section \ref{exis}.
Take $u_1-u_2$ as test function in each equation \eqref{resolv} (that  with datum $f_1$ and that with datum $f_2$) and subtract them. Then we get
\begin{eqnarray*}
&&\int_\Omega(u_1-u_2)^2dx+c\alpha\int_\Omega(u_1-u_2)^2dx\\
&+&c\int_\Omega(\mathbf z_1-\mathbf z_2,D(u_1-u_2))=\int_\Omega(f_1-f_2)(u_1-u_2) \,dx\,.
\end{eqnarray*}
Dropping a nonnegative term and applying H\"older's inequality, it follows that
\[ (1+c \alpha)\int_\Omega(u_1-u_2)^2dx\le\|f_1-f_2\|_2\|u_1-u_2\|_2\,, \]
from where we conclude
\[ \|u_1-u_2\|_2\le \frac{1}{1+c\alpha}\|f_1-f_2\|_2\,. \]
Therefore, if $\alpha>0$, then the Lipschitz constant satisfies $\frac{1}{1+c\alpha}<1$ and so each resolvent is actually a contraction.
Similar, simpler arguments show that the same result is true for
$1<p< 2$:
%\begin{multline}\label{Gp}
%G_p (u)=\\
%\left\{
%\begin{array}{ll}
%\displaystyle
%\displaystyle \frac1p \int_\Omega |\nabla u |^{p}\, dx + \alpha \int_\Omega u^2 dx\,,&\hbox{if }u\in W^{1,p}(\Omega)\cap L^2(\Omega)\,;\\[5mm]
%+\infty\,,&\hbox{if }u\in L^2(\Omega)\backslash W^{1,p}(\Omega)\,.
%\end{array}
%\right.\nonumber
%\end{multline}
\begin{equation}\label{Gp}
G_p (u)=
\left\{
\begin{array}{ll}
\displaystyle
\displaystyle \|\nabla u \|^{p}_p + \alpha \| u\|_2^2\,,&\hbox{if }u\in W^{1,p}(\Omega)\cap L^2(\Omega)\,;\\[5mm]
+\infty\,,&\hbox{if }u\in L^2(\Omega)\backslash W^{1,p}(\Omega)\,.
\end{array}
\right.
\end{equation}
\section{Numerical Resolution}
In this section we shall exploit the underlying structure of the minimization problem to write the corresponding energy functional as the difference of convex functions.
For this we consider functionals (\ref{minFun}) and (\ref{ep}) defined as $E_1 (u)=J_1 (u)+ H(u,f)$ and $ E_p (u)=J_p (u)+ H(u,f)$. Using (\ref{Guno}) and (\ref{Gp}) we can
decompose them in form $E_1 (u)=G_1 (u)- F(u,f)$ and $ E_p (u)=G_p (u)- F(u,f)$ where (compare with (\ref{hener}))

\begin{equation}\label{fener}
F (u,f)=
\lambda \displaystyle \int_\Omega \log I_0 \left( \frac{u f}{\sigma^2} \right)dx\,.
\end{equation}
The fundamental point is that the energy in (\ref{fener}) is convex.
As a consequence (\ref{minFun}) and (\ref{ep}) are difference of convex energy functionals.

We now introduce a 2D discrete setting in which the functionals can be minimized by a convergent Proximal Point algorithm,
in which a primal-dual method is used to solve the proximal operator for (\ref{Guno}) and (\ref{Gp}) together with an {\em ascent} gradient step for (\ref{fener}).
The generalization to 3D (volumetric) data sets is straightforward.

\subsection{Discrete Framework}
Let $\Omega \subset \mathbb{R}^2$ be an {\em ideally} continuous
rectangular image domain and consider a discretization in
terms of a regular Cartesian grid $\Omega_h$ of size $N \times M$:
$(ih, jh), \,1 \le i \le N,\, 1\le j \le M$ where $h$ denotes the
size of the spacing. The matrix $(u^h_{i,j})$ represents a
discrete image where each pixel $u_{i,j}$ is located in the
correspondent node $(ih, jh)$. In what follows, we shall choose $h
= 1$ because it only causes a rescaling of the energy through the
$\lambda$ parameter. Henceforth we shall drop the dependence of
the mesh size and denote $u^h =u$. Let $X=\mathbb{R}^{N \times M}$
be the space of solutions. We introduce the discrete gradient
$\nabla\,:\,X\to Y=X \times X$, defined as the forward finite differences operator:
\begin{equation}\label{disc_grad}
(\nabla u)_{i,j} =
\displaystyle\begin{pmatrix}\; (\nabla u)_{i.j}^x\;\\ \;(\nabla u)_{i.j}^y\;\end{pmatrix} =
\displaystyle\begin{pmatrix} \;u_{i+1, j} - u_{i, j}\;\\ \; u_{i, j+1} - u_{i, j}\;\end{pmatrix}
\end{equation}
except for $(\nabla u)^x_{N, j}=0$, and $(\nabla u)^y_{i, M} =0$. The discrete  p-norm of the gradient for $1\leq p <2$ is:
\[ \|\nabla u\|_p^p=\sum_{i.j} |(\nabla u)_{i.j}|^p \mbox{, with }\, \] \[ |(\nabla u)_{i.j}| = \sqrt{((\nabla u)_{i.j}^x)^2+((\nabla u)_{i.j}^y)^2} \]
which for $p=1$ is the discrete version of the isotropic TV operator (\ref{tv}) and for $1< p <2$, is the discrete version of the $J_p(u)$ term of the energy (\ref{ep}).
The discrete energy for the functionals defined in (\ref{minFun}) and (\ref{ep}) reads as:
\begin{equation}\label{ener_disc}
E_p(u)=\displaystyle \frac{1}{p} \sum_{i.j} |(\nabla u)_{i.j}|^p+ \frac{\lambda}{2\sigma^2}\sum_{i. j} u_{i, j}^2 - \lambda\sum_{i. j}\log I_0\left(\frac{u_{i,j}  f_{i, j}}{\sigma^2}\right) ,\;
\end{equation}
where the matrix $(f_{i, j})$ represents the discrete noisy image, with each pixel $f_{i, j}$ located at the node $(i, j)$.

% In the same way we can define the discrete version of the ROF problem deduced in (\ref{enrof}) which is (at a generic time step $t^n$)

% \begin{equation}\label{ener_discROF}
%  \sum_{i.j} |(\nabla u)_{i,j}|+ \left(\frac{1}{2\beta}\right) \sum_{i. j}\left( u_{i, j} - g^n_{i, j}\right)^2
% \end{equation}
% The algorithm presented in \cite{Chambolle04} is based in the dual formulation of the ROF problem then
Endowing the spaces $X$ and $Y$ with the standard Euclidean scalar product, the adjoint operator of the discrete gradient
(\ref{disc_grad}) is $\nabla^*=-\mbox{div}$.
%\begin{equation}\label{divp}
%<\nabla u, p>_Y = - <u, \mbox{div}\, p>_X
%\end{equation}
Given $p=(p^{x},\, p^{y}) \in Y$, we have
\[ (\mbox{div}\, p)_{i,j}=(p_{i, j}^{x}-p_{i-1, j}^{x})+(p_{i, j}^{y}-p_{i, j-1}^{y}) \]
for $2\leq i, j \leq N-1$. The term ($p_{i, j}^{x}-p_{i-1, j}^{x}$) is replaced with $p_{i,j}^{x}$ if $i=1$ and with $-p_{i-1,j}^{x}$ if $i=N$, while the term ($p_{i, j}^{y}-p_{i, j-1}^{y}$) is replaced with $p_{i, j}^{y}$ if $j=1$ and with $-p_{i, j-1}^{y}$ if $j=N$.
%Let moreover $\Gamma_0$ be the set of convex proper l.s.c functions
\subsection{A Proximal Point Algorithm for Rician Denoising} %Rician Denoising functional as a DC function: a Proximal Point Algorithm}
In this section we address the numerical resolution of the non-smooth non-convex minimization problem associated to the energy functional (\ref{minFun}) ($p=1$)
and the smooth non-convex approximating minimization problems related to the differentiable energy (\ref{ep}) ($1<p<2$).
To this end, we shall adapt a general proximal point algorithm for the minimization of the difference of  convex (DC) functions
proposed in \citep{Sun2003}. A decomposition of the energy functional as a difference of  convex (DC) functions is then proposed.
This is based on the fact that $I_0(s)$ is strictly log-convex which means that
$\log I_0 (s) $ is strictly convex and so is the energy term defined in (\ref{fener}).

In the discrete setting introduced before we can then write the Rician denoising functional (\ref{ener_disc}) as follows.
Given $f$, let $F : X\to \mathbb{R}$
and $G_p : X\to \mathbb{R}$
be the discretized analogue of functionals (\ref{fener}) and (\ref{Guno}) ($p=1$), (\ref{Gp}) ($1<p<2$):
\[ F(u)=\lambda \displaystyle \sum_{i. j} \log I_0\left(\frac{u_{i,j}  f_{i, j}}{\sigma^2}\right) , \mbox{ and}\]
\[ G_p (u)
%=\displaystyle \frac{1}{p}\|\nabla u\|_p^p+\frac{\lambda}{2\sigma^2}\|u\|_{\ell^2}^2
%=\frac{1}{p} \|\nabla u\|_p^p+\displaystyle \frac{\lambda}{2\sigma^2}\|u\|_{2}^2
=\displaystyle \frac{1}{p} \sum_{i.j} |(\nabla u)_{i.j}|^p+ \frac{\lambda}{2\sigma^2}\sum_{i. j} u_{i, j}^2 ,\quad 1\leq p < 2 \]
The functional in (\ref{ener_disc}) can be seen as the difference of two strictly convex proper l.s.c
functions $G_p(u)$ and $F(u)$:
\[ E_p (u)=G_p (u)-F(u) ,\qquad 1\leq p < 2 \]
Notice that $G_p(u) \geq 0$, $F(u) \geq 0$ and $F$ is differentiable with {\em Frechet derivative} $F' (u)$.

%Notice that $P_k =(I+c_k
%\,\partial G_p)^{-1} $ is a proximal mapping and define a Proximal
%Point type iteration
%$$
%u_{k+1}=P_k (u_k +c_k w_k ).
%$$
Then we can find a global minimizer of $E_p(u)$ by applying the following Proximal Point algorithm:
% Meter Box Algorithm
\begin{itemize}
	\item Given an initial  point $u_0=f$, let $c_k= c, \forall k$ and set  $k=0$ and $\epsilon=10^{-6}$.
	\begin{enumerate}\label{algPP}
		\item Compute $w_k= F' (u_k)$.
		\item Set $y_k=u_k+c_k w_k$ .
		\item Compute $u_{k+1}=(I+c_k \,\partial G_p)^{-1} (y_k)$
		\item If $\|u_{k+1} -u_k\|_2/\|u_k\|_2<\epsilon$ stop. Otherwise $k=k+1$ and return to step 1.
	\end{enumerate}
\end{itemize}
Notice that we can write Steps 1-3 as
\begin{align*}
u_{k+1}&=(I+c_k \,\partial G_p)^{-1} (u_k + c_k F'(u_k))\\
&=\mbox{Prox}_{c_k G_p}(u_k + c_k F'(u_k))\\
&=\mbox{Prox}_{c_k G_p}(u_k - c_k \partial (-F(u_k)))
\end{align*}
which is a forward-backward splitting algorithm (see for example \citep{Zhang2010a}).

Step 1 is explicitly given by
\[ w_k=F' (u_k) = \displaystyle \left(\frac{\lambda }{\sigma^2 }\right)%\frac{I_1(u_k f/\sigma^2)}{I_0(u_k f/\sigma^2)}f(x)
\displaystyle \frac{I_1 \displaystyle \left( \frac{u_k f }{\sigma^2 }\right) }{I_0 \displaystyle \left( \frac{u_k f}{ \sigma^2 } \right)} f \]
In Step 2 we set the {\em descent} direction for Step 3. Notice
that any {\em ascent} direction for $F$ is a {\em descent}
direction for $E$. To compute  the proximal operator $(I+c
\,\partial G_p)^{-1} $ in Step 3 we need to solve the following
strictly convex minimization problem:
\begin{eqnarray}\label{eq:ROFminNorm}
u_{k+1}&=\displaystyle \arg \min_{u\in X}& \displaystyle \left(
G_p(u) + \displaystyle \frac{1}{2c} \|u-y_k\|_{2}^2 \right)=\\
&=\displaystyle \arg\min_{u\in X} &\displaystyle \left(
%\frac{1}{p} \|\nabla u\|_p^p+\displaystyle \frac{\lambda}{2\sigma^2}\|u\|_{2}^2 + \frac{1}{2c} \|u-y_k\|_{2}^2\right)
\displaystyle \frac{1}{p} \sum_{i.j} |(\nabla u)_{i.j}|^p+\displaystyle \frac{\lambda}{2\sigma^2}\sum_{i. j} u_{i, j}^2\right.\nonumber\\
&&\left.+\frac{1}{2c}\sum_{i. j} (u_{i, j}-y^k_{i,j})^2\right) %\|u-y_k\|_{2}^2\right)\nonumber
\end{eqnarray}
%=\min_{u\in X} \|\nabla u\|_{\ell^1}+\frac{\lambda}{2\sigma^2}\|u\|_{\ell^2}^2+ \frac{1}{2c} \|u-y_k\|_{\ell^2}^2
%The non differentiability of $G_1 (u)$ is dealt with the Primal
%Dual algorithm presented in \cite{Chambolle2011}.
Let $R_p (u)=\displaystyle  \frac{1}{p} \sum_{i.j} |(\nabla u)_{i.j}|^p %\frac{1}{p}\|\nabla u\|_p^p , \, 1<p<2 ,
$, and \[S(u)=
\displaystyle  \frac{\lambda}{2\sigma^2}\sum_{i. j} u_{i, j}^2+ \frac{1}{2c}\sum_{i. j} (u_{i, j}-y^k_{i,j})^2 %\frac{\lambda}{2\sigma^2}\|u\|_{2}^2+ \frac{1}{2c} \|u-y_k\|_{2}^2 %,\quad \forall p
\]
Using Legendre Fenchel's duality,
we write the minimization problem (\ref{eq:ROFminNorm}) as a saddle point problem:
\begin{eqnarray}\label{eq:saddle}
&\arg \min_{u\in X}& \,\displaystyle \left(R_p (u)+S(u)\right)\\ =&\arg\min_{u\in X}& \,\displaystyle \left(\max_{v\in Y} \,\langle \nabla u, v\rangle - R_{p}^*(v)\right) +S(u)\nonumber
\end{eqnarray}
% F=R G=S
We distinguish two cases. When $p=1$, the Fenchel conjugate $R_{p}^*(v)$
is the indicator function $I_K$ of the convex set $K=\{ v \in
Y\,:\, \|v\|_{\infty} \leq 1\}$, i.e $I_K (v)=0$ if $v\in K$, $I_K
(v)=+\infty$ if $v\notin K$. In the differentiable case,
$1<p<2$, we have
\[ R_{p'}^*(v)=\displaystyle \frac{1}{p'}\|v\|_{p'}^{p'} =\displaystyle  \frac{1}{p'}\sum_{i. j} |v_{i, j}|^{p'}  \]
with $1/p+1/p'=1$. To solve this saddle-point problem (\ref{eq:saddle}) we use the Primal Dual algorithm presented in \citep{Chambolle2011}. This method allows an unified treatment of (\ref{eq:saddle}) for any $p$, so dealing with the non differentiability of $G_1 (u)$.
This algorithm performs Step 3 in k$^{th}$ external iteration of the Proximal Point
algorithm and reads as follows:

Given $u^0=y_k$, set $v^0 =\bar{0}$, $\tau_d=\tau_p=1/\sqrt{12}$ and $\bar{u}^0=u^0$. Iterate until convergence:
\begin{enumerate}
	\item[(i)] $v^{n+1}=\left( I + \tau_d \partial R_{p'}^* \right)^{-1} (v^n + \tau_d \nabla \bar{u})$
	\item[(ii)] $u^{n+1}=(I + \tau_p \partial S)^{-1} (u^n + \tau_p \,\mbox{div}\, v^{n+1} )$
	\item[(iii)] $\bar{u}^{n+1}=2u^{n+1}-u^n$
\end{enumerate}
This is an inner loop and the upper index $n$ is the inner iteration counter. Steps (i) and (ii) aim to compute the proximal operators corresponding to $R_p^*(u)$ and $S(u)$ and are defined by:

Step (i). For $p=1$, we compute $\bar{v}^n =v^n+\tau_d \nabla \bar{u}^n$ and the resolvent operator with respect to $R_1^*$ reduces to pointwise Euclidean projector onto $\ell^2$ balls:
\[ v^{n+1}=\displaystyle  \left( I + \tau_d \partial R^* \right)^{-1} (\bar{v}^n ) \quad\Longleftrightarrow \quad v_{i, j}^{n+1} = \frac{\bar{v}_{i, j}^{n}}{\max(1, |\bar{v}_{i,j}^{n}|)} \]
For $1<p<2$, with $\bar{v}^n =v^n+\tau_d \nabla \bar{u}^n$ the computation of the resolvent operator $\left( I + \tau_d \partial R^* \right)^{-1} (\bar{v}^n ) $ leads to solve  the following strictly convex minimization problem:
\begin{eqnarray*}
v^{n+1}&=& \displaystyle \arg \min_{v\in Y} \displaystyle \left(
\frac{1}{p'}\|v\|_{p'}^{p'}+\frac{1}{2\tau_d} \|v-\bar{v}^n
\|_{2}^2 \right)\\
&=& \displaystyle \arg \min_{v\in Y} \displaystyle \left( \frac{1}{p'}\sum_{i. j} |v_{i, j}|^{p'} +
\frac{1}{2\tau_d }\sum_{i. j} (v_{i, j}-\bar{v}_{i,j})^2 \right)
\end{eqnarray*}
The first order necessary (and sufficient) condition for optimality reads:
\[ f^n (v)=\tau_d |v|^{p'-2}v+v-\bar{v}^n=0 \]
It is easily seen that $f^n (v)$ is continuous, monotone increasing with $f^n (0)=\bar{v}^n$ and the
equation $f^n (v)=0$ has a unique real positive solution $0<| v| \leq |\bar{v}^n|$ for any $p'$.
For any fixed internal iteration $n$ we apply the Newton's method to solve the nonlinear equation resulting in the following fixed point iteration:
Set $j=0$, $v_{j}^{k,n+1} =v_{j}$, $v_{j+1}^{k,n+1} =v_{j+1}$ and $v_0 =v^{k,n}$. Compute, for $j=1, 2, ...$ till convergence
\[ v_{j+1}=\phi (v_j ) =%v_k-\displaystyle\frac{\tau_d |v_k|^{p'-2}v_k+v_k-\bar{v} }{ \tau_d (p'-1)|v_k |^{p'-2}+1} =
\displaystyle\frac{\tau_d (p' -2)|v_j|^{p'-2}v_j+\bar{v} }{ \tau_d (p'-1)|v_j |^{p'-2}+1} \]

Step (ii). The resolvent operator with respect to $S$ poses simple
pointwise quadratic problems. The solution is  given by
\[ u=\displaystyle (I + \tau_p \partial S)^{-1} (\bar{u} ) \quad \Longleftrightarrow \quad u_{i,j}=%\frac{y_k/c + \bar{u}/\tau_p}{(\lambda/\sigma^2) + c^{-1}+\tau_p^{-1}}
\frac{\sigma^2 \left(\tau_p y_k + c\bar{u}\right)_{i,j}}{c \tau_p \lambda + \sigma^2 (c+\tau_p)} \]
%In Step (iii) we actualize the intermediate variable $\bar{u}$ and
%then we get back to Step 4 in the Proximal Point algorithm to set
%a new {\em descent} direction until convergence.
\section{Numerical Results}
In this section we test the performance of the proposed numerical scheme. We first validate the results of the TV Rician denoising method using the Proximal Point Algorithm (PPA), denoted by TV-Rician in the following. We also test the numerical convergence of the $p$-approximating problems. Then, we compare TV-Rician with previously proposed methods for TV Rician-based denoising \citep{Martin2011,Getreuer2011,Chen2015} for different images and noise intensities. Finally we present an application on real Diffusion Tensor Images (DTI), which is an MRI modality heavily affected by Rician noise \citep{Basu2006,Tristan-Vega2010}.
\subsection{Numerical Scheme Validation}
\begin{figure*}[h!tb]
	\begin{center}
		{\includegraphics[width=0.72\textwidth]{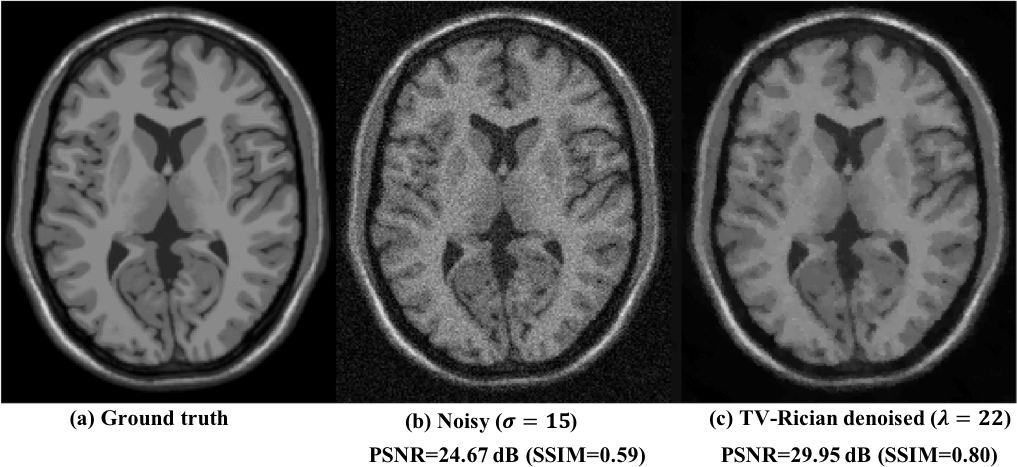}}	
	\end{center}
	\caption{Denoising test on a phantom brain image. At left, the original free-of-noise slice. In the center, the same slice contaminated with Rician noise for $\sigma=15$. At right, the best denoised image obtained using TV-Rician as measured by PSNR and SSIM.}\label{im:Riciandenoising}
\end{figure*}

\begin{figure*}[h!tb]
	\begin{center}
		{\includegraphics[width=\textwidth]{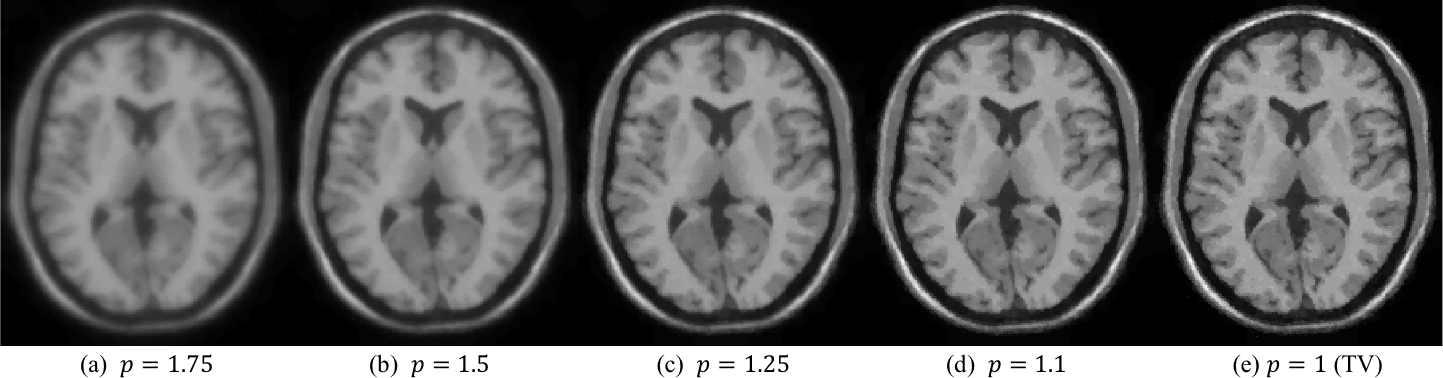}}	
	\end{center}
	\caption{Denoising results of the noisy phantom brain image of Figure \ref{im:Riciandenoising} using the p-Laplacian for $p=1.75, 1.5, 1.5, 1.1$ and $1$ (Total Variation).}\label{im:Lapdenoising}
\end{figure*}

In order to assess the performance of the proposed algorithm we used a synthetic brain image obtained from the BrainWeb Simulated Brain Database \footnote{available at
	http://www.bic.mni.mcgill.ca/brainweb} at the Montreal Neurological Institute \citep{Cocosco1997}.
The central slice of the original phantom was extracted and normalized to be between 0 and 255. Finally, the slice was contaminated artificially with Rician noise for $\sigma=15$. To compute the denoising quality we use two different measures: the Peak-Signal-to-Noise-Ratio (PSNR) and the Structural Similarity Index (SSIM) \citep{Wang2004}.

\begin{figure}[h!tb]
	\begin{center}
		{\includegraphics[width=0.75\columnwidth]{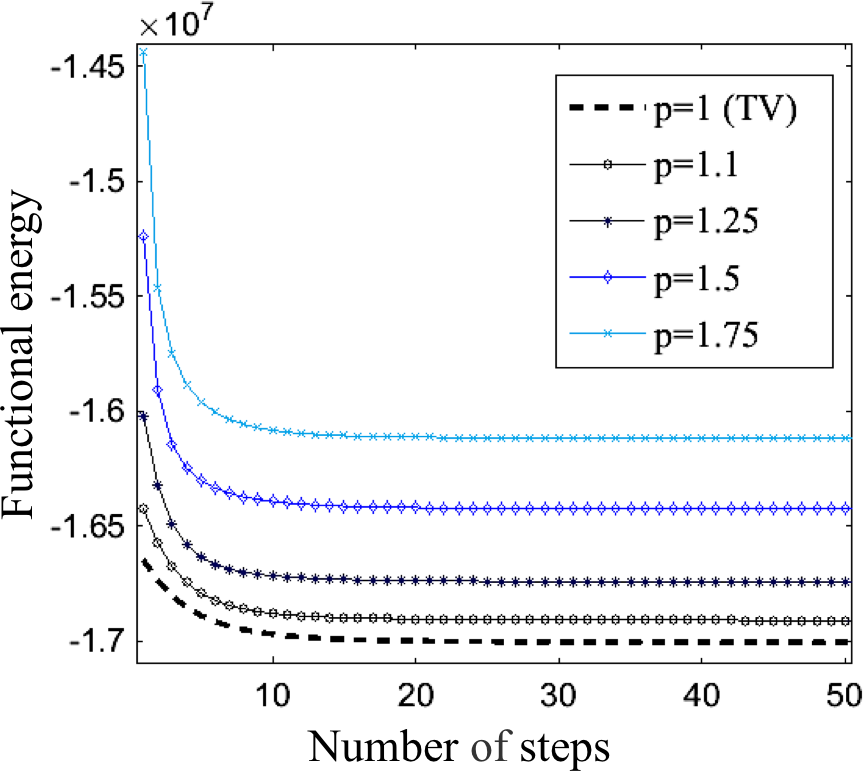}}
	\end{center}
	\caption{Energy minimization evolution of functional \eqref{ener_disc} for $p=1.75, 1.5, 1.5, 1.1$ and $1$ (Total Variation) for the images displayed in Fig. \ref{im:Lapdenoising}.}\label{im:energiesPlap}
\end{figure}

In Figure \ref{im:Riciandenoising}, we show the denoising results of the TV-Rician method for $\lambda=22$. \textcolor{black}{This $\lambda$ value was optimized to obtain the best PSNR and SSIM with respect to original phantom.} We can see how in the denoised image (Fig. \ref{im:Riciandenoising}c) most of the noise has been removed while the fine details are preserved. Using the same regularization parameter, we repeat this test solving \eqref{ener_disc} for different values of $p$, $p=\{1.1, 1.25, 1.5, 1.75\}$, to numerically asses the convergence of the $p-$sequence of regularizing approximating $u_p$ solutions when $p\to 1$. The $u_p$ solutions are shown in Figure \ref{im:Lapdenoising}, where, as expected, the  closer $p$ gets to 1, the more similar the $p$-Laplacian solution is to the TV image. This $p$-convergence can also be observed when plotting the energy minimization evolution of the Proximal Point Algorithm for these same values of $p$ and the TV case (see Figure \ref{im:energiesPlap}).
\subsection{Comparison with Other Variational Methods for Rician Denoising}
In order to cope with the difficulties of the non-smooth non-convex problem \eqref{minFun}, several methods have been proposed for TV-based denoising of Rician contaminated images. The first of them uses an $\epsilon$-approximation of the TV term \citep{Martin2011,Getreuer2011} to obtain a smooth minimization problem. With this regularization, a gradient descent can be applied to solve the problem. In the following this approach will be denoted as  {TV$_\epsilon$-Rician}. In the work of \citep{Getreuer2011}, a convexification of the functional was also proposed. This new minimization problem is solved by a Split-Bregman approach \citep{Goldstein2009a}. We will refer to this convexification as {Getreuer} model in the following. Finally, a different convexification of \eqref{minFun} by adding the term $\frac{1}{\sigma}\int_{\Omega}(\sqrt{u}-\sqrt{f})^2 dx)$ has been recently presented in \citep{Chen2015}. This new convex problem is then efficiently solved using a primal-dual algorithm \citep{Chambolle2011}. For the comparisons, we will denote as {Chen-Zeng} this method.

All of these approaches rely on approximations of the problem \eqref{minFun} making it differentiable or convex. Notably, the proposed algorithm ({TV-Rician}) based on the PPA scheme copes with the original non-smooth non-convex functional. For this comparison we use four images kindly provided by the authors of \citep{Chen2015}: one natural image
Camera man (256$\times$256), and three MR images Lumbar-Spine (200$\times$200), Brain ($217\times181$) and Liver (214$\times$304). The images are then corrupted by Rician noise for $\sigma=20$ and $\sigma=30$. For the sake of fairness, all the algorithms were run until fulfill the same convergence criterium based on the relative difference between the functional energy in two consecutive iterations. In our test, we set the tolerance to $1\times10^{-7}$. TV$_\epsilon$-Rician, Getreuer and TV-Rician use a regularization parameter $\lambda$ which multiplies the data fidelity term, while the Chen-Zeng algorithm uses a parameter $\gamma=1/\lambda$ multiplying the TV term. For all tests, the regularization parameters were separately optimized to get the best PSNR and to get the best SSIM with respect to the original images. The results of this comparison are displayed in Table \ref{table:results}. Notice that the optimal regularization parameter for each case is displayed in the table as $\gamma$.
\begin{table*}[h!tb]
	\caption{Comparisons of the best PSNR values and SSIM values by different methods for Rician denoising}\label{table:results}
	\begin{tabular}{llllll}
		\hline  &  &$\sigma=20$  &  &$\sigma=30$  &  \\
		\cline{3-6}Image  &Method  &PSNR $(\gamma=1/\lambda)$  &SSIM $(\gamma=1/\lambda)$   &PSNR $(\gamma=1/\lambda)$   &SSIM $(\gamma=1/\lambda)$  \\
		\hline Camera man  &TV$_\epsilon$-Rician  &28.12 (0.03)  &0.8077 (0.05)  &24.81 (0.02)  &0.7148 (0.025)  \\
		&Getreuer  &27.83 (0.03)  &0.7478 (0.05)  &25.58 (0.02)  &0.6653 (0.03)  \\
		&Chen-Zeng  &28.44 (0.035)  &0.8229 (0.045)  &25.69 (0.025)  &0.7539 (0.035)  \\
		&TV-Rician  &\textbf{28.64} (0.03)  &\textbf{0.8272} (0.04)  &\textbf{26.18} (0.025)  &\textbf{0.7655} (0.03)  \\
		
		\hline Lumbar-Spine  &TV$_\epsilon$-Rician  &28.27 (0.03) &0.7716 (0.03)  &25.28 (0.015)  &0.6609 (0.02)  \\
		&Getreuer  &27.66 (0.035) &0.6685 (0.035)  &24.81 (0.02)  &0.5115 (0.02)  \\
		&Chen-Zeng  &28.35 (0.035) &0.7743 (0.04)  &25.53 (0.025)  &0.6705 (0.03)  \\
		&TV-Rician  &\textbf{28.84} (0.03) &\textbf{0.7892} (0.053)  &\textbf{26.5} (0.02)  &\textbf{0.6998} (0.02)  \\
		
		\hline Liver  &TV$_\epsilon$-Rician  &29.06 (0.03)  &0.8033 (0.03)  &26.61 (0.02)  &0.7201 (0.025)  \\
		&Getreuer  &28.84 (0.035)  &0.7723 (0.04)  &26.75 (0.025)  &0.6742 (0.025)  \\
		&Chen-Zeng  &29.25 (0.04)  &0.8047 (0.04)  &27.03 (0.03)  &0.7371 (0.03)  \\
		&TV-Rician  &\textbf{29.4} (0.035)  &\textbf{0.8088} (0.035)  &\textbf{27.42} (0.025)  &\textbf{0.7452} (0.025)  \\
		
		\hline Brain  &TV$_\epsilon$-Rician  &26.7 (0.025)  &0.6550 (0.035)  &23.67 (0.015)  &0.5881 (0.025)  \\
		&Getreuer  &\textbf{29.41} (0.035)  &\textbf{0.8996} (0.04)  &\textbf{27.03} (0.025)  &\textbf{0.8000} (0.03)  \\
		&Chen-Zeng  &26.63 (0.035)  &0.6634 (0.04)  &23.79 (0.025)  &0.6067 (0.03)  \\
		&TV-Rician  &28.12 (0.03)  &0.6780 (0.04)  &25.61 (0.02)  &0.6165 (0.025)  \\
		
		\hline Average  &TV$_\epsilon$-Rician  &28.04 (-)  &0.7594 (-)  &25.09 (-)  &0.6710 (-)  \\
		&Getreuer  &28.43 (-)  &0.7721 (-)  &26.04 (-)  &0.6628 (-)  \\
		&Chen-Zeng  &28.17 (-)  &0.7663 (-)  &25.51 (-)  &0.6921 (-)  \\
		&TV-Rician  &\textbf{28.75} (-)  &\textbf{0.7758} (-)  &\textbf{26.42} (-)  &\textbf{0.7068} (-)  \\
	\end{tabular}
\end{table*}

\begin{figure*}[h!tb]
	\begin{center}
		\subfloat[I (original image)]{\includegraphics[width=0.19\textwidth]{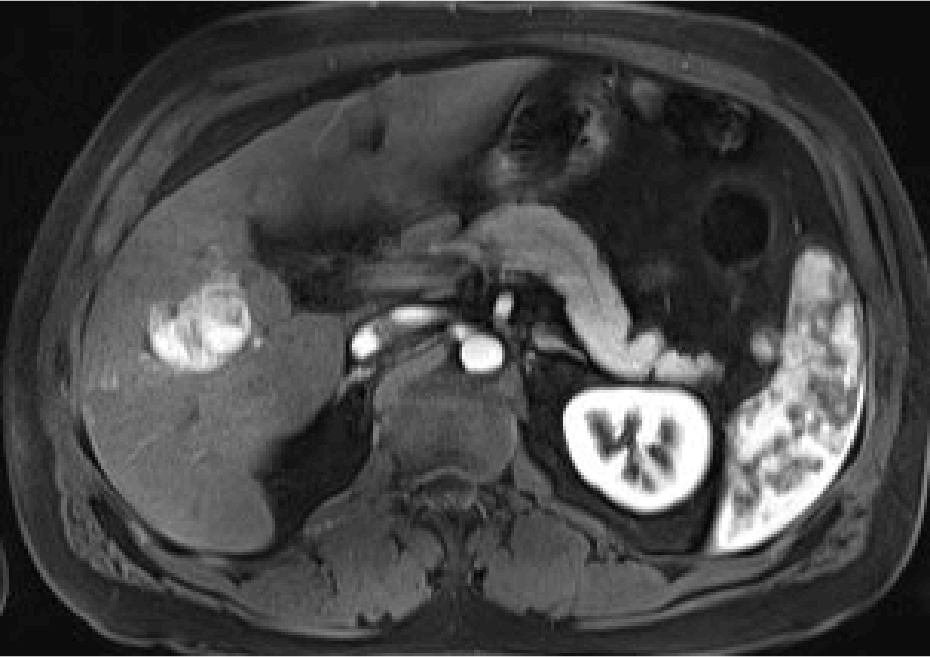}}\hfill
		\subfloat[TV$_\epsilon$-Rician]{\includegraphics[width=0.19\textwidth]{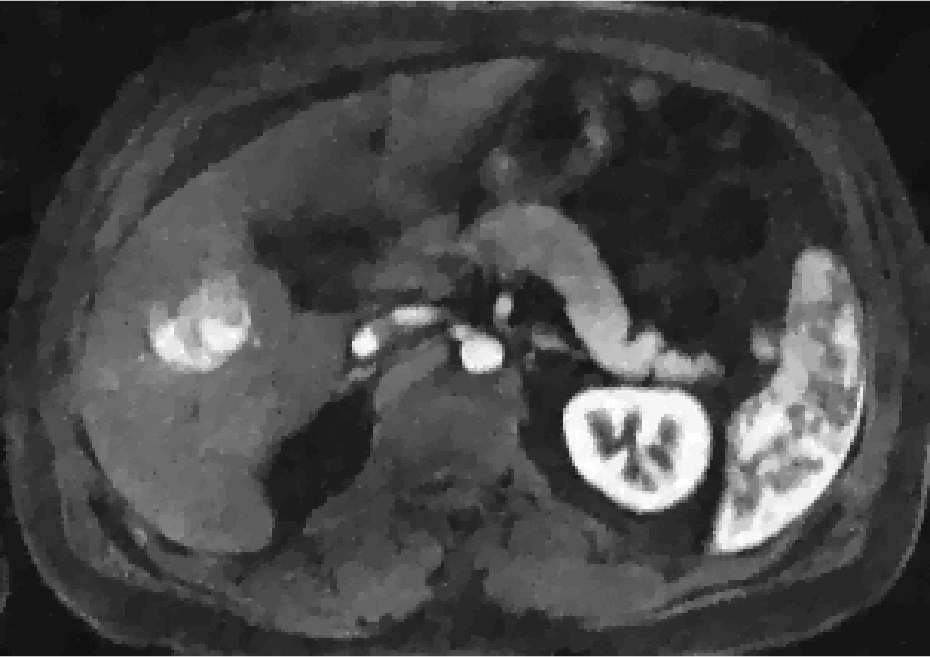}}\hfill
		\subfloat[Getreuer]{\includegraphics[width=0.19\textwidth]{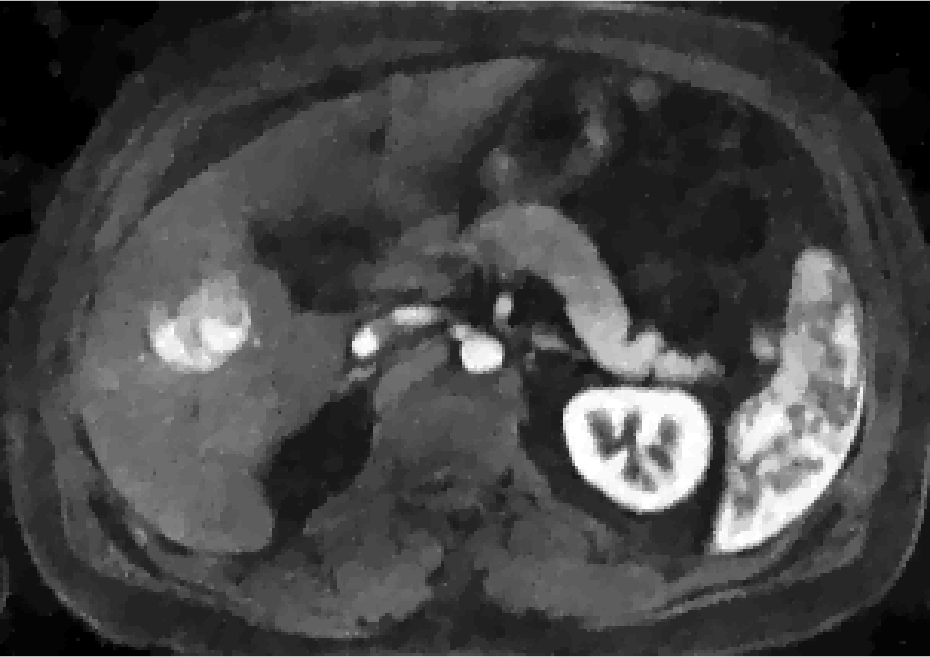}}\hfill
		\subfloat[Chen-Zang]{\includegraphics[width=0.19\textwidth]{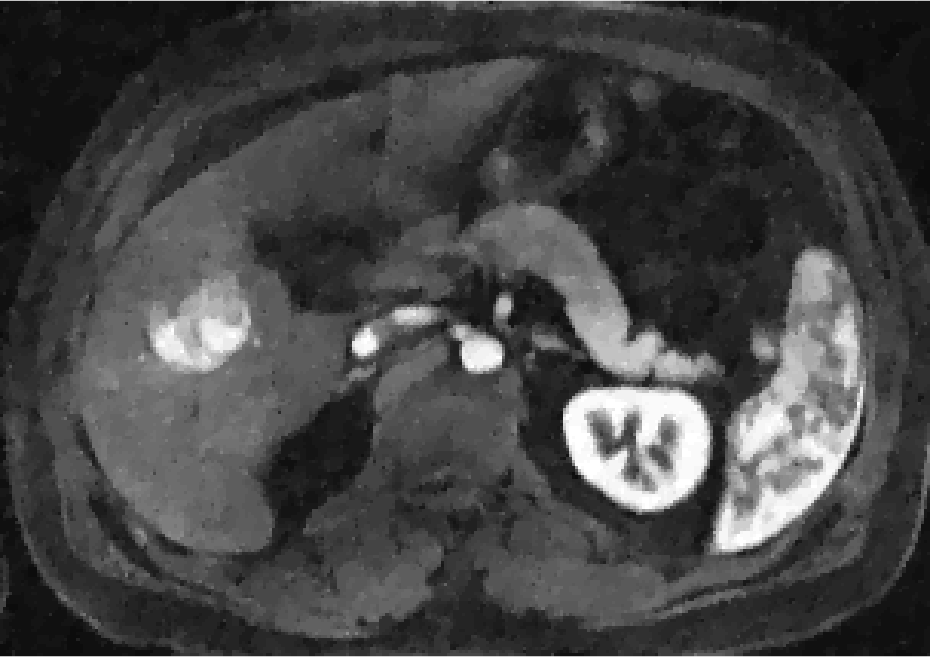}}\hfill
		\subfloat[TV-Rician]{\includegraphics[width=0.19\textwidth]{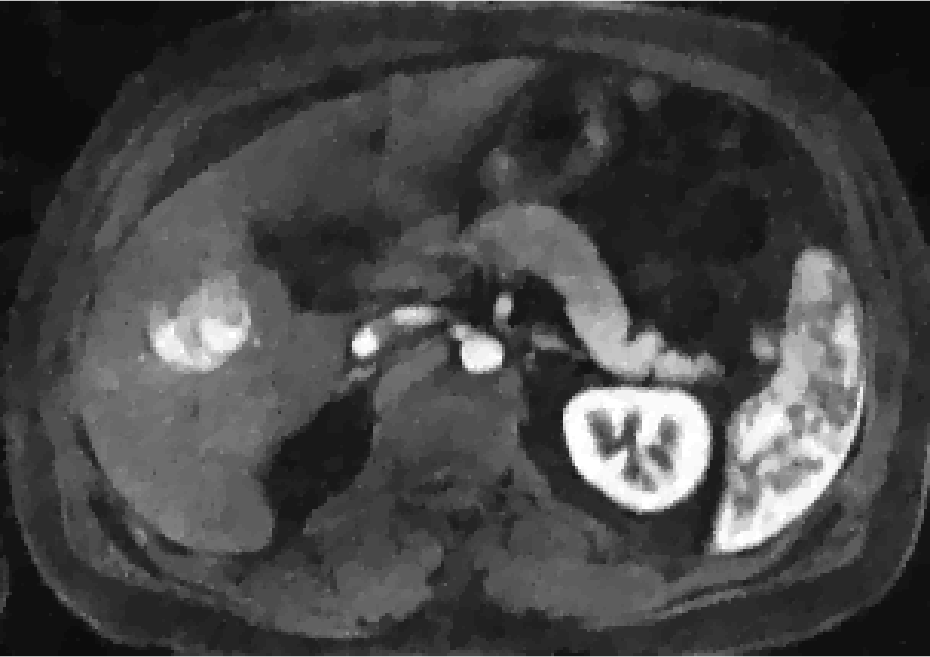}}\hfill
		\subfloat[f (noisy with $\sigma=20$) ]{\includegraphics[width=0.19\textwidth]{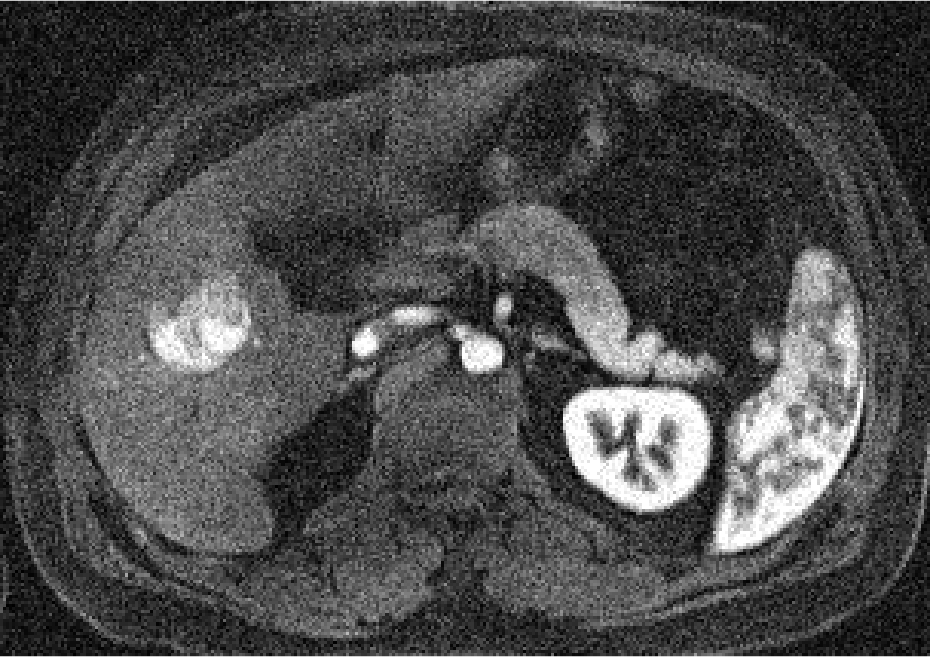}}\hfill
		\subfloat[$|$ TV$_\epsilon$-Rician - I $|$]{\includegraphics[width=0.19\textwidth]{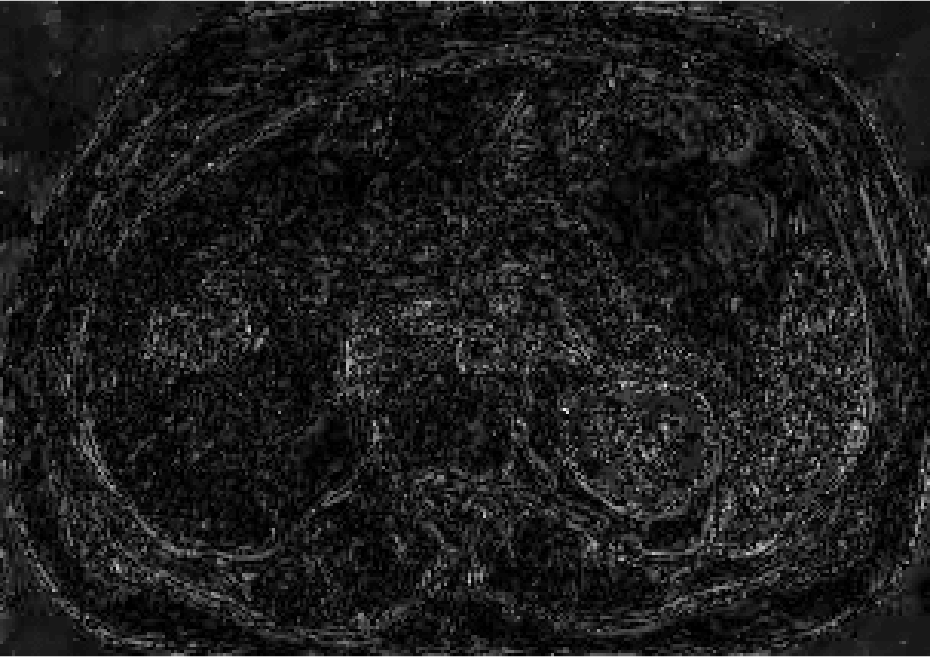}}\hfill
		\subfloat[$|$ Getreuer - I $|$ ]{\includegraphics[width=0.19\textwidth]{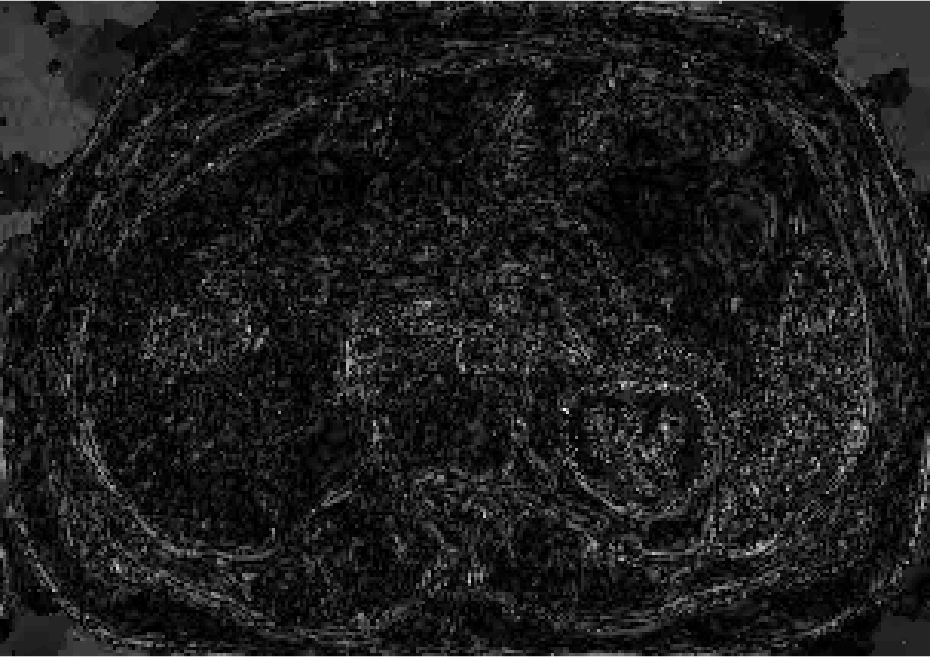}}\hfill
		\subfloat[$|$ Chen-Zang - I $|$]{\includegraphics[width=0.19\textwidth]{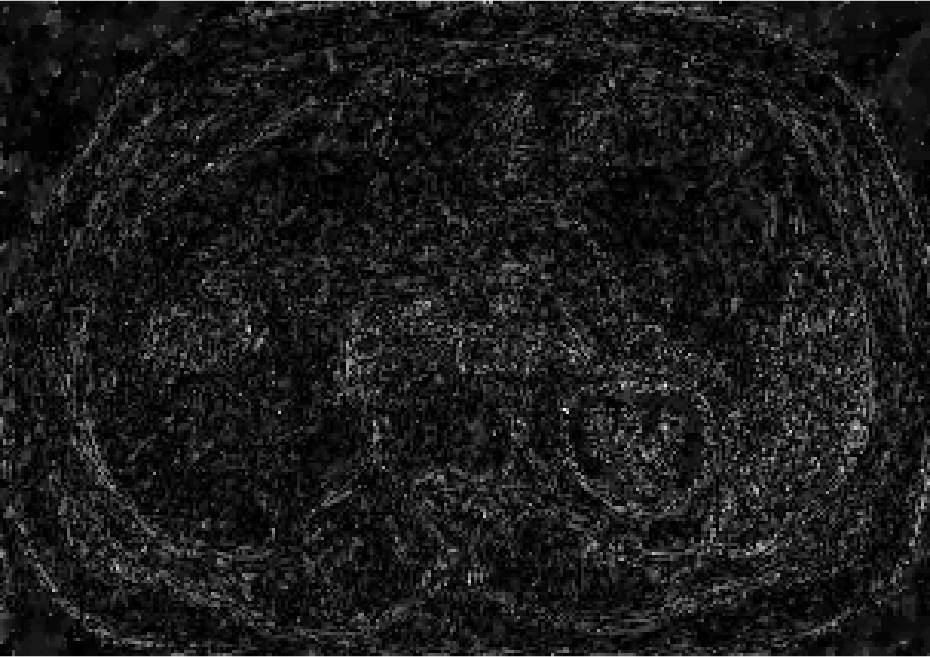}}\hfill
		\subfloat[$|$ TV-Rician - I $|$ ]{\includegraphics[width=0.19\textwidth]{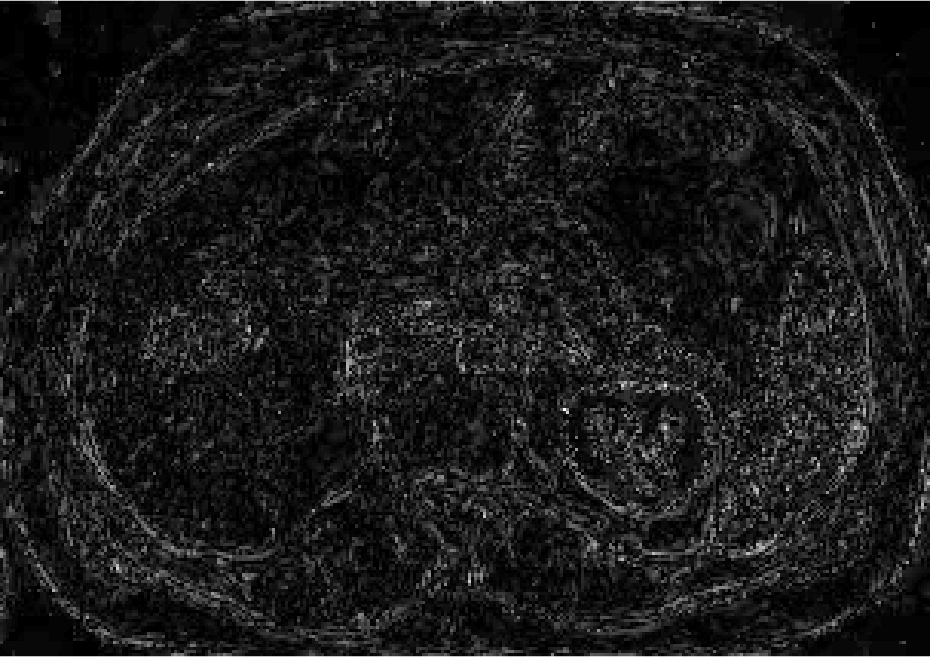}}\hfill
	\caption{Denoising results on the Liver image (a) for fixed parameters $\sigma=20$, $\gamma=0.035$.
		In (b)--(e), the images resulting from applying the compared methods to $f$ (f), the noisy version of $I$ (a), the original image. In the second row, (g)--(j), the absolute differences between the denoised images and $I$ are shown. Careful inspection reveals a better perfomance of the proposed method (see (e) and (j)) in areas with lower SNR (dark zones).}
		\label{im:comparisonRic}
	\end{center}
\end{figure*}

We see that TV-Rician  gets the best results in both PSNR and SSIM for the Camera man, the Lumbar-Spine and the Liver images for all levels of noise. The differences with other methods increase for higher noise level ($\sigma=30$), confirming that the original problem \eqref{minFun} is best suited than its approximations for Rician denoising. \textcolor{black}{For the case of the phantom Brain image, Getreuer model scores the best denoising results. In order to convexify the Rician data fidelity term the authors in \citep{Getreuer2011} substitute the original functional for small values of the solution by a linear approximation. This modified functional drives these values of the solution closer to 0 than the original Rician functional we considered. Since the background of the synthetic Brain image is 0, this model achieves a better solution for this image than the other methods. This effect can be observed in the other images. For instance, in Fig. \ref{im:comparisonRic}h, the error of this model in the upper corners is considerably higher than in the rest of the algorithms because the background in the noise-free Liver image is not 0.} Nevertheless, the proposed method (TV-Rician) gets higher PSNR and SSIM than TV$_\epsilon$-Rician and Chen-Zeng in the synthetic Brain image, and it is the best algorithm overall when computing the averaged PSNR and SSIM.

 Moreover, when using the same regularization parameter for all the methods, the TV-Rician method also achieves a solution which is a lower minimum of \eqref{minFun}. This comparison is performed for the Liver image and the parameters $\sigma=20$ and $\gamma=0.035$. These results are shown in Figure \ref{im:comparisonRic} and Table \ref{table:energ}: TV-Rician achieves the best denoising solution in terms of visual inspection, PSNR and energy minimization.
\begin{table}[h!tb]
	\caption{PSNR and energy functional values (see \eqref{minFun}) for the compared methods for Rician denoising.}\label{table:energ}
	\begin{tabular}{llll}
		\hline  &  &$\sigma=20, \gamma=0.035$  & \\
		\cline{3-4}
		Image  &Method  &PSNR  &$E_1$ \eqref{minFun}\\
		\hline Liver  &TV$_\epsilon$-Rician  &28.98 &$-1.4183\times10^7$  \\
		&Getreuer  &28.84  &$-1.4184\times10^7$\\
		&Chen-Zeng  &29.18 &$-1.4180\times10^7$\\
		&TV-Rician  &29.4 &$-1.4192\times10^7$\\
\end{tabular}
\end{table}
	
\subsection{Application on Real Diffusion Tensor Imaging of the Brain}	
\begin{figure*}[h!tb]
	\begin{center}
		{\includegraphics[width=\textwidth]{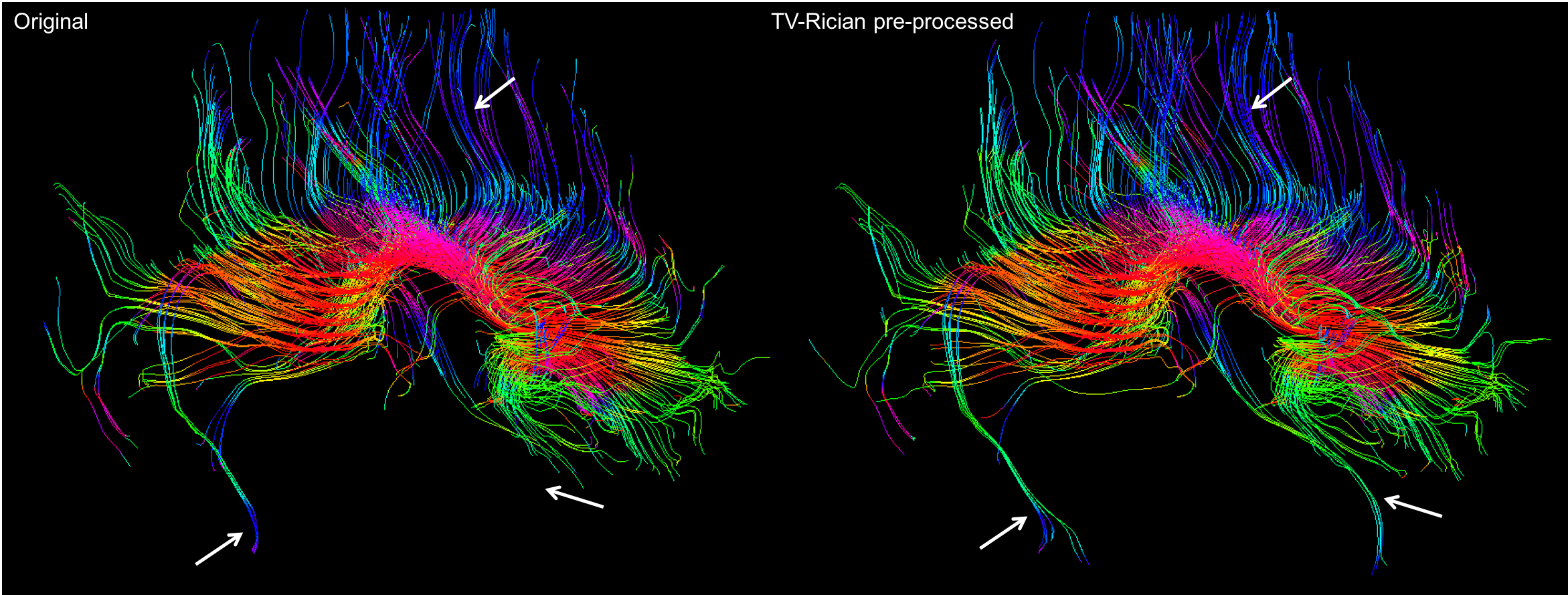}}
	\end{center}
	\caption{Tractography generated from a seed placed in the corpus callosum. At left, the tractography generated from the original DWI (and DTI) data. At right, the tractography generated from the TV-Rician denoised data.  Particular areas where the tractographies are different because of the noise are pointed by white arrows in the image}\label{im:wholeTract}
\end{figure*}

\begin{figure*}[h!tb]
	\begin{center}
		{\includegraphics[width=\textwidth]{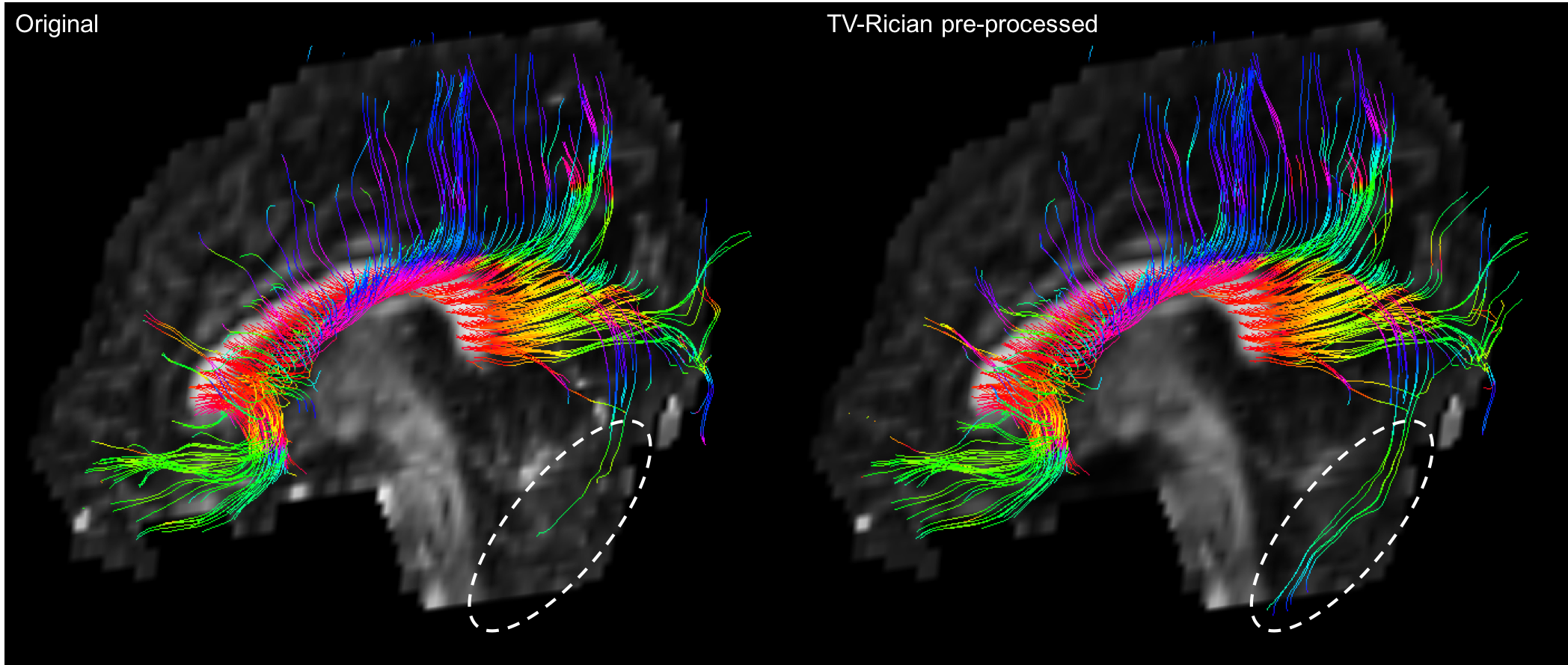}}
	\end{center}
	\caption{Tractography generated from the corpus callosum over a sagital view of the FA image. Notice that after the pre-processing, the tracts of the left arcuate fasciculus is recovered.}\label{im:tract_Arcuate}
\end{figure*}

The data we used consist of a Diffusion Weighted Images (DWI) dataset provided by Fundaci\'on CIEN-Fundaci\'on Reina Sof\'ia which was acquired with a 3 Tesla General Electric scanner equipped with an 8-channel coil. The DWI have been obtained with a single-shot spin-echo EPI sequence (FOV = 24 cm, TR = 9600 ms, TE = 91.5 ms, slice thickness = 2 mm, spacing = 0.6 mm, matrix size = 128x128, NEX = 1). The DWI data consists on a volume obtained with b=0 s/mm$^2$ and 45 volumes with b=1000 s/mm$^2$ corresponding with gradient directions that equally divide the 3-D space. These DWI, which represent diffusion measurements along multiples directions, are denoised by solving the proposed minimization problem \eqref{minFun} using the PPA.
Then, Diffusion Tensor Images (DTI) are reconstructed from the original and denoised DWI data using the 3D Slicer tools\footnote{Freely available in http://www.slicer.org/}. DTI is one of the most popular methods for in vivo analysis of the white matter (WM) structure of the brain, helping to detect WM alterations that can be found from early stages in some degenerative diseases \citep{Gattellaro2009}. The DTI information is commonly used to generate a tractography of a particular area of the brain, which is a 3D representation of the fibers of WM involved. In Figure \ref{im:wholeTract}, the tractographies generated from a seed placed in the corpus callosum are shown. White arrows indicate regions where the noise in the image generated from the original data (at left) affects to the reconstruction of specific tracts which are nevertheless recovered in the tractography from the pre-processed data (at right).
In order to highlight the regions where the fibers reconstruction differs we display the tractographies over a sagital view of the Fractional Anistropy (FA) generated from the same DTI data (Fig. \ref{im:tract_Arcuate}). It can be seen how the left arcuate fasciculus can not be reconstructed from the original data but it is recovered after the pre-processing. The correct reconstruction of the left arcuate fasciculus is important since it is involved in important tasks like language and praxis \citep{Catani2008}.

\section{Conclusions}
In this paper we presented the mathematical analysis of the quasi-linear elliptic equation for the $1$-laplacian operator which arises from considering the minimization of the Total Variation based energy functional modeling Rician denoising for MRI. Theoretical difficulties come from both ingredients of the model: the TV regularization term, which makes the problem non--smooth, and the Rician statistics of the noise in MRI, which yields a non-convex minimization problem.
We provided sufficient conditions on the data for the existence of a bounded non-trivial BV solution of the elliptic equation which turns out to be a global minimizer of the associated energy functional. Several qualitative properties of this solution have been deduced.
The uniqueness of a strictly positive solution is still an open problem. Extensive numerical experiments not reported here suggest that there exists only one such solution.

We also proposed and implemented a convergent Proximal Point Algorithm to solve this non--smooth non--convex minimization problem.  The numerical results demonstrate the effectiveness of the proposed method compared to previous approximations to TV-based Rician denoising. Finally, we tested our algorithm in in-vivo DTI tractography showing the benefits of pre-processing  DWI data before DTI reconstruction.
%\begin{acknowledgements}
%\end{acknowledgements}



\bibliographystyle{spbasic}
\bibliography{jmiv_MSS_2015}

\end{document}